\RequirePackage{fix-cm}
\documentclass[smallextended]{svjour3}       
\smartqed  
\usepackage{graphicx}
\usepackage{latexsym}
\usepackage{calc}
\usepackage{datetime}
\usepackage{amssymb,amsmath}

\usepackage{graphicx}

\usepackage{url}
\usepackage{hyperref}

\usepackage{subcaption}
\usepackage[misc,geometry]{ifsym} 

\usepackage{booktabs}

\usepackage{xcolor}

\newcounter{hours}\newcounter{minutes}

\def\Def{\stackrel{\mathrm{def}}{=}}

\def\inter{{\rm int \,}}

\def\dom{{\rm dom \,}}
\def\beq{\begin{equation}}
\def\eeq{\end{equation}}

\def\C{\mathbb{C}}
\def\E{\mathbb{E}}
\def\Q{\mathbb{Q}}
\def\R{\mathbb{R}}

\def\e{\varepsilon}

\newcommand{\refLE}[1]{\ensuremath{\stackrel{(\ref{#1})}{\leq}}}
\newcommand{\refEQ}[1]{\ensuremath{\stackrel{(\ref{#1})}{=}}}
\newcommand{\refGE}[1]{\ensuremath{\stackrel{(\ref{#1})}{\geq}}}

\newcommand{\half}{\mbox{${1 \over 2}$}}

\def\ba{\begin{array}}
\def\ea{\end{array}}
\def\beann{\begin{eqnarray*}}
\def\eeann{\end{eqnarray*}}
\def\bea{\begin{eqnarray}}
\def\eea{\end{eqnarray}}

\def\BT{\begin{theorem}}
\def\ET{\end{theorem}}
\def\BL{\begin{lemma}}
\def\EL{\end{lemma}}
\def\BC{\begin{corollary}}
\def\EC{\end{corollary}}
\def\BE{\begin{example}}
\def\EE{\end{example}}
\def\BD{\begin{definition}}
\def\ED{\end{definition}}
\def\BR{\begin{remark}}
\def\ER{\end{remark}}
\def\BAS{\begin{assumption}}
\def\EAS{\end{assumption}}
\def\BI{\begin{itemize}}
\def\EI{\end{itemize}}

\def\BMP{\begin{minipage}{9.5cm}}
\def\EMP{\end{minipage}}
\def\MPT{\begin{minipage}{11.5cm}}
\def\EPT{\end{minipage}}

\def\la{\langle}
\def\ra{\rangle}

\begin{document}

\title{Improved global performance guarantees of second-order methods in convex minimization
}

\titlerunning{Improved global complexity of second-order methods in convex minimization}        

\author{Pavel Dvurechensky{$^\dagger$}\thanks{$\dagger$Correspondence to  \email{pavel.dvurechensky@wias-berlin.de}}        \and
        Yurii Nesterov 
}


\institute{P. Dvurechensky \at
              Weierstrass Institute for Applied Analysis and Stochastics, Mohrenstr. 39, 10117 Berlin, Germany \\
            \url{orcid.org/0000-0003-1201-2343}\\  
			\email{pavel.dvurechensky@wias-berlin.de}           
           \and
           Yu. Nesterov \at
              Center for Operations Research and Econometrics (CORE), Catholic University of Louvain (UCL), 34 voie du Roman Pays, 1348 Louvain-la-Neuve, Belgium\\
              \email{Yurii.Nesterov@uclouvain.be} 
}

\date{Received: date / Accepted: date}

\maketitle

\begin{abstract}
In this paper, we attempt to compare two distinct branches of research on second-order optimization methods. The first one studies self-concordant functions and barriers, the main assumption being that the third derivative of the objective is bounded by the second derivative. The second branch studies cubic regularized Newton methods (CRNMs) with the main assumption that the second derivative is Lipschitz continuous. We develop a new theoretical analysis for a path-following scheme (PFS) for general self-concordant functions, as opposed to the classical path-following scheme developed for self-concordant barriers. We show that the complexity bound for this scheme is better than that of the Damped Newton Method (DNM) and show that our method has global superlinear convergence. 
We propose also a new predictor-corrector path-following scheme (PCPFS) that leads to further improvement of constant factors in the complexity guarantees for minimizing general self-concordant functions.
We also apply path-following schemes to different classes of constrained optimization problems and obtain the resulting complexity bounds.
Finally, we analyze an important subclass of general self-concordant functions, namely a class of strongly convex functions with Lipschitz continuous second derivative, and show that for this subclass CRNMs give even better complexity bounds.
\end{abstract}

\section{Introduction}
\label{intro}
\vspace{1ex}\noindent
{\bf Motivation.} Local performance guarantees for the
second-order methods are known since 1948 \cite{kantorovich1948newton}. In that paper, the author proved a local
quadratic convergence of the Newton method under some
natural assumptions (non-degeneracy of the Hessian at
solution and local Lipschitz continuity of the Hessian).
However, in some sense, the quadratic convergence is too
fast: each step of such methods doubles the number of
right digits in the approximate solution. Therefore,
questions on the acceleration of these schemes were never
raised in the literature \cite{conn2000trust}. Moreover, for many years the only global complexity results for second-order methods were
obtained in the framework of the theory of self-concordant
functions and barriers \cite{nesterov1994interior,nesterov2004introduction}.

The situation was changed after the paper \cite{nesterov2006cubic}, where
the first global complexity bounds were obtained for the cubic regularized Newton method (CRNM). Namely,
it was shown that for convex functions with globally
Lipschitz continuous Hessian the CRNM converges in terms of the
function value as $O(\frac{1}{k^2})$, where $k$ is the
iteration counter. Very soon it was shown that this method
can be accelerated up to the rate $O(\frac{1}{k^3})$ using
the technique of estimating sequences  \cite{nesterov2008accelerating}. Under the same assumptions, based on a special line-search procedure, the authors of \cite{monteiro2013accelerated} proposed a second-order method with convergence rate $O\left({\frac{1}{k^{7/2}}}\right)$.

Thus, at this moment there exist two main, nearly independent, frameworks for global complexity analysis of the second-order methods. One is based on the affine-invariant theory of self-concordant functions. And the second one assumes bounded third derivatives of the objective in a fixed Euclidean norm. In the recent years, many new results were obtained within the second framework, including development of accelerated schemes.
The main goal of this paper is to revisit the first framework and try to improve existing complexity bounds. The secondary goal of this paper is to show that these two classes of problems do intersect and we can compare the efficiency of the corresponding methods. Towards the first goal, we derive new complexity bounds for a path-following scheme (PFS) as applied to the unconstrained minimization of a self-concordant {\em function}. This result is new since the known complexity bounds for path-following methods are related to self-concordant {\em barriers} (see, for example, Section 4.2 in \cite{nesterov2004introduction}). Moreover, this result shows that our PFS exhibits {\em global superlinear convergence}, and the resulting complexity bound is better than that of the classical Damped Newton Method (DNM). 
We next propose a predictor-corrector path-following scheme (PCPFS) that leads to further improvement of constant factors in the complexity guarantees for minimizing self-concordant functions. Additionally, we show that our results lead to improved complexity for constrained problems where one can additionally use the self-concordant barrier property.
We further compare our bounds with the complexity results for different versions of the CRNM on the class of strongly convex functions with Lipschitz continuous Hessian. It appears that such functions are self-concordant. We conclude that the latter methods are much more efficient when applied to this  subclass of self-concordant functions.

\vspace{1ex}\noindent
{\bf Related work.} Local convergence rate analysis of Newton method was extended in \cite{lee2014proximal} to the case of composite optimization with the objective given as a sum of a twice-differentiable function and a simple convex function. Inexact proximal Newton methods for self-concordant functions with composite term  were studied in \cite{li2017inexact}, where global and local convergence were studied. For the Machine Learning applications, the authors of \cite{rodomanov2016superlinearly} propose an  incremental Newton method with local superlinear convergence and global linear convergence. The authors of \cite{zhang2018communication} study distributed minimization of self-concordant empirical losses via DNM combined with distributed preconditioned conjugate gradient method.  A new version of Newton method with adaptive stepsizes was proposed in \cite{polyak2019new} for non-linear systems of equations. 
Generalization of DNM for generalized self-concordant functions is developed together with local convergence analysis in \cite{sun2018generalized}. 
The authors of \cite{hanzely2022damped} analyze global and local convergence of a DNM for minimizing a subclass of self-concordant functions called semi-strongly self-concordant functions.
Beyond classical second-order Newton-type methods, local superlinear convergence was obtained for path-following interior-point methods \cite{nesterov2016local}. 
The papers \cite{dinh2013inexact,liu2020inexact} study inexact path-following method for optimization with special separable structure.
In \cite{tran-dinh2018new}, the authors propose a new homotopy proximal variable-metric framework for composite convex minimization and show that under appropriate assumptions such as strong convexity-type and smoothness, or self-concordance, their schemes can achieve global linear convergence. 
A general Adaptive Regularization algorithm using Cubics was proposed in \cite{cartis2011adaptive}, where global and local convergence were proved under relaxed assumptions. Global complexity analysis was extended for the case of H\"older-continuous Hessians in \cite{grapiglia2017regularized} and even higher order derivatives in \cite{cartis2019universal,grapiglia2020tensor,song2021unified}. Linear convergence of the CRNM under H\"older continuity of the Hessians and uniform convexity was proved in \cite{doikov2021minimizing}. The author of \cite{doikov2023minimizing} obtains global complexity bounds for minimizing quasi-self-concordant functions by gradient regularization of Newton method.
Frank-Wolfe method for minimization of self-concordant functions was recently proposed in \cite{staudigl2020self-concordant}. A combination of Frank-Wolfe and Newton method for minimization of self-concordant functions was proposed in \cite{liu2020newton}, and close approach, but for functions with Lipschitz-continuous Hessian, was studied in \cite{carderera2020second-order}. In \cite{dvurechensky2024hessian} the authors combine the ideas of cubic regularization and self-concordant barriers to propose algorithms for constrained non-convex optimization with global complexity guarantees.

Previous works analyze either the class of self-concordant functions or the class of functions with Lipschitz-continuous Hessians. In this work, we attempt to compare these two classes in terms of the global complexity of solving the corresponding minimization problem. At the same time, previous works analyzed the path-following methods for self-concordant barriers rather than for self-concordant functions. Here we do not assume the barrier property and analyze  path-following schemes for general self-concordant functions. In particular, to the best of our knowledge, this is the first time when global superlinear convergence is obtained for this class of methods and problems.

{\bf Contributions.} Our contributions can be summarized as follows.
\begin{enumerate}
	\item We give a new theoretical analysis of a path-following scheme (PFS) for \textit{general} self-concordant \textit{functions}. This is in contrast to \cite{nesterov2004introduction}, where the classical analysis is made under additional assumption that the objective is self-concordant \textit{barrier}. Moreover, the obtained global complexity bound of the PFS is better than that of the Damped Newton Method (DNM). In particular, we demonstrate global superlinear convergence of the PFS.
	\item We consider a feasibility problem and compare complexity bounds for this problem by the DNM, PFS and dual PFS. Importantly, the proposed PFS gives the best complexity, improving upon known complexities.
	\item We propose a predictor-corrector path-following scheme (PCPFS) for \textit{general} self-concordant \textit{functions}. This scheme has even better global complexity bound in terms of the constant factors than the PFS for minimizing general self-concordant functions and also possesses superlinear convergence.
	\item We apply the PCPFS to two types of constrained optimization problems and under additional self-concordant barrier property show that the resulting primal and dual PCPFS have better constant factors in the complexity compared to the existing bounds in \cite{nesterov2004introduction}.
	\item We propose variants of our path-following schemes with adaptive stepsizes that may take longer steps adapting to local behavior of the central path.
	\item We make an observation that the above two classes of functions intersect. Namely, strongly convex functions with Lipschitz continuous Hessian belong to both classes. Thus, all the discussed methods can be applied to this class. Our conclusion is that on this class CRNMs possess better complexity that the DNM and the PFS. 
\end{enumerate}

\vspace{1ex}\noindent
{\bf Contents.} In Section \ref{sc-SCF} we recall the properties of self-concordant functions extending them to the case of non-standard self-concordant functions and give the complexity of the DNM. In Section \ref{sc-PFS}, we provide description and new analysis of the PFS for general self-concordant functions. In Section \ref{S:feasibility} we discuss implications of our analysis and algorithms when applied to a feasibility problem. Section \ref{sc-PCPFS} is devoted to the construction and analysis  of our PCPFS for minimizing general self-concordant functions. In Section \ref{sc-PCPFS-appl}, we apply the PCPFS to two types of constrained optimization problems that involve self-concordant barriers and obtain the resulting complexity.
Section \ref{sc-FStrong} contains complexity analysis of strongly convex functions with Lipschitz-continuous Hessians by CRNMs and their comparison to the complexity by the DNM and PFS. 

\vspace{1ex}\noindent
{\bf Notation.}
Let $\E$ be a finite-dimensional real vector space and $\E^*$ be its dual.
For a function $f$, we denote its domain as ${\rm dom} f = \{x \in \E: f(x) < + \infty \}$, and the closure of the domain as ${\rm Dom } f$.
Given a function $f$ with  Hessian $\nabla^2f(x)$ that is non-degenerate at any $x \in \E$, we denote
$$
\ba{rcl}
\| h \|_x & = & \la \nabla^2 f(x) h,h \ra^{1/2}, \; h \in
\E, \quad \| g \|^*_x \; = \; \la g, [\nabla^2f(x)]^{-1}g
\ra^{1/2}, \; g \in \E^*
\ea
$$
\beq\label{def-Lambda}
\ba{rcl}
\text{and} \;\; \lambda_f(x) & = & \| \nabla f(x) \|^*_x, \quad x \in \E.
\ea
\eeq
We use $D^3f(x)[h]=\lim_{\alpha \to 0} \frac{1}{\alpha} (\nabla^2 f(x+\alpha h)-\nabla^2 f(x))$ to denote the third-order derivative at a point $x \in \E$ in a direction $h \in \E$.
For $h,h_1,h_2,h_3 \in \E$ we use a shortcut notation
\[
D^3f(x)[h]^3= \la D^3f(x)[h] h,h\ra, \quad D^3f(x)[h_1,h_2,h_3]= \la D^3f(x)[h_1] h_2,h_3\ra.
\] 
We define also, for $\tau \geq 0$, $\omega(\tau) = \tau - \ln(1 + \tau)$, $\omega_*(\tau) = - \tau - \ln(1-\tau)$.

\section{Minimizing self-concordant functions: Damped Newton Method}\label{sc-SCF}

We consider the following minimization problem:
\beq\label{prob-SCF}
\ba{c}
f^* \; = \; \min\limits_{x \in \E} f(x),
\ea
\eeq
where $f$ is a general self-concordant function whose definition is given below. We assume
that the Hessian of this function at any point is positive
definite and that the solution $x^*$ of problem
(\ref{prob-SCF}) exists.

\subsection{Preliminaries on general self-concordant functions}
In this technical subsection we give the definition and properties of general self-concordant functions which we use throughout the paper.
Let us start from a variant of the definition of
self-concordant functions.
\BD
Let function $f$ from $\C^3$, i.e., three times continuously differentiable, be convex on $\E$. It is
called a general self-concordant function if there exists
a constant $M_f \geq 0$ such that for any point $x \in \E$
and direction $h \in \E$:
\beq\label{def-SCF}
\ba{rcl}
|D^3f(x)[h]^3| & \leq & 2 M_f \la \nabla^2 f(x) h, h
\ra^{3/2}.
\ea
\eeq
If $M_f=1$, then the function is called {\em standard
self-concordant}. \qed
\ED
An equivalent characterization of general self-concordant functions is given by the following result \cite[Lemma 4.1.2]{nesterov2004introduction}.
\BL\label{lm-Eqiv-Def-SCF}
A function $f$ from $\C^3$ is general self-concordant iff for any point $x \in \E$
and any $h_1,h_2,h_3 \in E$ we have
\beq\label{def-SCF-2}
|D^3f(x)[h_1,h_2,h_3] | \leq  2 M_f \|h_1\|_x \|h_2\|_x \|h_3\|_x.
\eeq
\EL

It is clear that for any self-concordant function $f$,
function
\beq\label{eq-Norm}
\ba{rcl}
\tilde f (x) & = & M_f^2 f(x), \quad x \in \E,
\ea
\eeq
is standard self-concordant. Standard self-concordant
functions are more convenient for defining
self-concordant barriers (see \cite{nesterov1994interior}). However, for the main algorithms in this paper we work directly with definition (\ref{def-SCF}) and do not impose an additional assumption that $f$ is a self-concordant barrier.

Taking into account normalization (\ref{eq-Norm}), we can
rewrite all known properties of standard self-concordant
functions for the general ones. Let us present the most
important of them (see Section 4.1 in \cite{nesterov2004introduction}).
For all $y \in \E$ with $\| y - x \|_x < {1 \over M_f}$ we
have
\beq\label{eq-Compat}
\ba{rcl}
(1 - M_f \| y - x \|_x)^2 \nabla^2 f(x) & \preceq &
\nabla^2 f(y) \; \preceq \; {1 \over (1 - M_f \| y - x
\|_x)^2} \nabla^2 f(x),
\ea
\eeq
\beq\label{eq-Norm-Compat}
\ba{rcl}
 \| y - x \|_y & \leq & \frac{\|y-x\|_x}{1-M_f\|y-x\|_x},
\ea
\eeq
\beq\label{eq-UBound}
\ba{rcl}
f(y) & \leq & f(x) + \la \nabla f(x), y - x \ra + {1 \over
M_f^2} \omega_*(M_f \| y - x \|_x),
\ea
\eeq
and for all
$y \in \E$ we have
\beq\label{eq-LBound}
\ba{rcl}
f(y) & \geq & f(x) + \la \nabla f(x), y - x \ra + {1 \over
M_f^2} \omega(M_f \| y - x \|_x).
\ea
\eeq
Inequality (\ref{eq-LBound}) leads to the following bound for $\lambda_f(x) < {1 \over M_f}$
\beq\label{eq-FBound}
\ba{rcl}
f(x) - \min\limits_{y \in \E} f(y) & \leq {1 \over M_f^2} \omega_*(M_f \lambda_f(x)).
\ea
\eeq
Moreover, for $\lambda_f(x) < {1 \over M_f}$ and $x^*=\arg \min \limits_{y \in \E} f(y)$,
\beq\label{eq-x-lambda-Bound}
\ba{rcl}
\|x-x^*\|_x & \leq {1 \over M_f} \omega_*'(M_f \lambda_f(x)) = {M_f \lambda_f(x) \over 1-M_f \lambda_f(x)}.
\ea
\eeq
Similarly to \eqref{eq-Compat}, if
$\delta \equiv \| \nabla f(x) - \nabla f(y) \|^*_x < {1
\over M_f}$, then
\beq\label{eq-CompD}
\ba{rcl}
(1 - M_f \delta)^2 \nabla^2 f(x) & \preceq & \nabla^2 f(y)
\; \preceq \; {1 \over (1 - M_f \delta)^2} \nabla^2 f(x).
\ea
\eeq
This relation follows from the observation \cite{nesterov1994interior} that
the Fenchel conjugate $f_*$ for $f$ defined as
$$
\ba{rcl}
f_*(s) & = & \sup\limits_{x \in \E} [ \la s, x \ra -
f(x)]
\ea
$$
is also self-concordant on its domain with the same
constant $M_f$. Note also that, for $x(s) = \arg \max \limits_{x \in \dom f} [ \la s, x \ra -
f(x)]$,
\begin{align}
&\nabla f_*(s) = x(s),  \label{eq-Grad-Dual} \\
&\nabla^2 f_*(s) = [\nabla^2 f(x(s))]^{-1}. \label{eq-Hess-Dual}
\end{align}

The following result was not proved in \cite{nesterov2004introduction} even for the case of $M_f=1$. 
\BL\label{lm-scf-grad-diff}
Let $x,y \in {\rm dom} f$ be such that $\| y - x \|_x < {1 \over M_f}$. Then
\[
\|\nabla f(y)-\nabla f(x) \|_x^* \leq {\|y-x\|_x \over 1- M_f\|y-x\|_x}.
\]
\EL
\proof
The proof follows from inequality (4.1.7) in \cite{nesterov2004introduction}, Cauchy-Schwarz and the duality relation \eqref{eq-Grad-Dual}.
\qed
%
%
%

\subsection{Newton methods for general self-concordant functions.}

We now move to standard algorithms used for minimizing general self-concordant functions.
When applied to this class of functions, Standard Newton Method has local quadratic convergence.
\BT[Theorem 4.1.14 in \cite{nesterov2004introduction}]\label{lm-StandNStep} Define the Standard Newton Step
\beq\label{eq-StandNStep}
\ba{rcl}
x_+ & = & x - [\nabla^2f(x)]^{-1} \nabla f(x).
\ea
\eeq
Then,
\beq\label{eq-StandQuad}
\ba{rcl}
\lambda_f(x_+) & \leq & {1 \over M_f} \left( {M_f \lambda_f(x)} \over {1- M_f \lambda_f(x)} \right)^2.
\ea
\eeq
\ET

Minimizing the right-hand side of inequality
(\ref{eq-UBound}) in $y$, we come to the following result.
\BT\label{lm-NStep} Define the Damped Newton Step
\beq\label{eq-NStep}
\ba{rcl}
x_+ & = & x - {[\nabla^2f(x)]^{-1} \nabla f(x) \over 1 +
M_f \lambda_f(x)}.
\ea
\eeq
Then,
\beq\label{eq-Decr}
\ba{rcl}
f(x_+) & \leq & f(x) - {1 \over M_f^2} \omega(M_f
\lambda_f(x)).
\ea
\eeq
Moreover,
\beq\label{eq-Quad}
\ba{rcl}
\lambda_f(x_+) & \leq & 2 M_f \lambda_f^2(x).
\ea
\eeq
\ET
\proof
Inequality (\ref{eq-Decr}) follows from Theorem 4.1.2 in \cite{nesterov2004introduction}.

We now prove inequality (\ref{eq-Quad}), which was not proved in \cite{nesterov2004introduction} even for $M_f=1$.
Denote $\lambda = \lambda_f(x)$, $h = x_+ - x$, $r = \| h \|_x =
{\lambda \over 1 + M_f \lambda}$. Hence, $\lambda = {r
\over 1 - M_f r}$. Note that, since $r < {1 \over M_f}$, 
$$
\ba{rcl}
\lambda^2_+ & \equiv & \la \nabla f(x_+), [\nabla^2
f(x_+)]^{-1} \nabla f(x_+) \ra \; \refLE{eq-Compat} \; {1
\over (1- M_f r)^2} \la \nabla f(x_+), [\nabla^2
f(x)]^{-1} \nabla f(x_+) \ra.
\ea
$$
Without changing the notation, we can associate with the
Hessians symmetric positive-definite matrices. Then,
denoting $G = [\nabla^2 f(x)]^{1/2} \succ 0$, we have
\beq
\label{eq-lm-NStep-proof-1}
\ba{rcl}
\nabla f(x_+) & = & \nabla f(x) + \int\limits_0^1 \nabla^2
f(x+\tau h) h d \tau \\
&\refEQ{eq-NStep} & -(1+M_f \lambda) \nabla^2 f(x) h +
\int\limits_0^1 \nabla^2 f(x+\tau h) h d \tau\\
& = & G \left[ -(1+M_f \lambda)I + G^{-1}
\left( \int\limits_0^1 \nabla^2 f(x+\tau h) d \tau \right) G^{-1}
\right] G h.
\ea
\eeq
By Corollary 4.1.4 in \cite{nesterov2004introduction}, we have that
\[
(1 - M_f r +
{1 \over 3} M_f^2 r^2) \nabla^2 f(x) \preceq \int\limits_0^1 \nabla^2 f(x+\tau h) d \tau \preceq {1
\over 1 - M_f r} \nabla^2 f(x).
\]
Thus, denoting $H = [ \cdot ]$ in \eqref{eq-lm-NStep-proof-1}, we can see that $H \preceq 0$ and 
$H \succeq [ -(1+M_f \lambda) + (1-M_f r)]I = [
-(1+M_f \lambda) + {1 \over 1 + M_f \lambda}]I \succeq - 2
\lambda M_f I$. So, we conclude that
$$
\ba{rcl}
\lambda_+^2 & \leq & {1 \over (1-M_f r)^2} (\| G H G h
\|^*_x)^2 \; = \; {1 \over (1-M_f r)^2} \la G H^2 G h, h
\ra \leq {(2 M_f \lambda )^2 \over (1-M_f r)^2} \la
G^2 h , h \ra \\
\\
& = & {(2 M_f \lambda )^2 \over (1-M_f r)^2}  \la
\nabla^2f(x) h, h \ra \; = \; {(2 M_f \lambda )^2 r^2\over
(1-M_f r)^2} \; = \; 4 M_f^2 \lambda^4.
\hfill{\mbox{\qed}}
\ea
$$

We can now analyze the efficiency of the Damped Newton Method (DNM)
\beq\label{met-DNM}
\ba{rcl}
x_{k+1} & = & x_k - { [\nabla^2 f(x_k)]^{-1} \nabla f(x_k)
\over 1 + M_f \lambda_f(x_k)}, \quad k \geq 0
\ea
\eeq
in terms of iteration complexity as applied to the minimization problem \eqref{prob-SCF}.
In view of inequality (\ref{eq-Quad}),
method (\ref{met-DNM}) starts converging quadratically when
it enters the region
\beq
\label{eq-quadr-conv}
\ba{rcl}
\Q & = & \left\{x \in \E: \; \lambda_f(x) \leq {1 \over
2M_f} \right\}.
\ea
\eeq
This convergence is very fast and, in view of inequality
(\ref{eq-FBound}), any reasonable accuracy in function
value can be reached in a small number of iterations.
Therefore, the main computational time is spent when
$\lambda_f(x_k) \geq {1 \over 2M_f}$. Denote by $N$ the last
iteration such that $\lambda_f(x_k)   \geq   {1 \over 2M_f}$,   $k = 0, \dots, N.$
Then, in view of inequality (\ref{eq-Decr}), we have
\beq\label{eq-BoundCN}
N \leq {\Delta(x_0) \over \omega\left(\half \right)} ,
\quad \Delta(x_0) \; \Def \; M_f^2 (f(x_0) - f^*).
\eeq
This is the main result of this section. Namely, when applied to the class of general self-concordant functions, DNM has complexity
$O(\Delta(x_0))$. In the following sections, we will show that this complexity may be improved by proposing a different algorithm.

Let us show that $\Delta(x_0)$ is a natural complexity
measure of our problem class. In order to see this, let us
attribute to our objects some physical units. Denote the
units for measuring the function value by $\mu_f$, and the
units for measuring the argument by $\mu_x$. Then, the
units for measuring the gradient are $\mu_g = \mu_f/\mu_x$.
The Hessian is measured in $\mu_h = \mu_f / \mu_x^2$, and
the third derivative is measured in $\mu_t = \mu_f/
\mu_x^3$. Thus, in view of definition (\ref{def-SCF}), the
units for measuring the constant $M_f$ are $\mu_s  =  \mu_t \mu_x^3/ (\mu_h \mu_x^2)^{3/2}  = 
\mu_f^{-1/2}.$
Note that the number of iterations is an integer number
with no physical dimension (scalar). Therefore, for using
the constant $M_f$ in the bounds for the number of
iterations, it must be multiplied by something having
physical dimension $\mu_f^{1/2}$. The simplest way to do
this is to define the characteristic $\Delta(x_0)$ as in
(\ref{eq-BoundCN}). In the sequel, we will use
$\Delta(x_0)$ as the main characteristic of complexity of
problem (\ref{prob-SCF}). Importantly, 
we can use the characteristics of our problem as arguments
of nonlinear univariate functions only by transforming
them to a scalar form. For example, the values $M_f
\lambda_f(x)$ and $M_f \| h \|_x$ have no physical
dimension and can be substituted to the functions $\omega$, $\omega_*$.

In addition to the global iteration complexity \eqref{eq-BoundCN}, we can show that the DNM accelerates to  superlinear convergence when $\lambda_f(x_k)\leq \frac{1}{M_f}$.
\BT\label{th-DNM-superlinear} 
For the DNM \eqref{met-DNM} when $\lambda_f(x_k)\leq \frac{1}{M_f}$ it holds that
\beq\label{eq-DNM-superlinear}
M_f^2(f(x_{k+1})-f^*) \leq M_f^2(f(x_{k})-f^*) \left( 1-\frac{\omega(\omega_*^{-1}(M_f^2(f(x_{k})-f^*)))}{M_f^2(f(x_{k})-f^*)} \right).
\eeq
In other words, the method has superlinear convergence.
\ET
\proof
Denote for $k\geq 0$ $\Delta_k=f(x_{k})-f^*$ and $\lambda_k=\lambda_f(x_k)$. Then, by \eqref{eq-FBound}, we have that $\omega_*^{-1}(M_f^2 \Delta_k) \leq  M_f \lambda_k$. Combining this with \eqref{eq-Decr}, we obtain 
\[
M_f^2 \Delta_{k+1} \leq M_f^2 \Delta_k - \omega(\omega_*^{-1}(M_f^2 \Delta_k)),
\]
which is \eqref{eq-DNM-superlinear}. Denoting $\omega_*^{-1}(M_f^2 \Delta_k) = \xi$, we see that $\frac{\omega(\omega_*^{-1}(M_f^2 \Delta_k))}{M_f^2 \Delta_k}=\frac{\omega(\xi)}{\omega_*(\xi)}$. Recalling the definition of $\omega(\xi), \omega_*(\xi)$ and using their Taylor expansion, we readily see that $\frac{\omega(\xi)}{\omega_*(\xi)}< 1$ and, hence, \eqref{eq-DNM-superlinear} leads to convergence. Further, since $\omega_*^{-1}$ monotonically increases, we have that $\xi$ decreases as $\Delta_k$ decreases to $0$. Considering the function $\frac{\omega(\xi)}{\omega_*(\xi)}$, its derivative, and Taylor expansions for $\omega(\xi)$, $\omega_*(\xi)$, it is easy to show that this function increases when $\xi$ decreases and its limit when $\xi\to 0$ is equal to $1$. Combining all the observations, we see that indeed the DNM has superlinear convergence when $\lambda_f(x_k)\leq \frac{1}{M_f}$, i.e., local superlinear convergence.
In the next section we propose a path-following scheme that has \textit{global superlinear convergence}.

\section{Minimizing self-concordant functions: path-following scheme}\label{sc-PFS}

In this section, we estimate the complexity of solving problem (\ref{prob-SCF}) by a path-following scheme (PFS). 
Note that the full justification of path-following methods is done so far in \cite{nesterov2004introduction}
only for self-concordant \textit{barriers}. On the contrary, our analysis relies only on the less restrictive assumption that $f$ is a self-concordant \textit{function}. Moreover, we make the analysis for \textit{general} self-concordant functions with $M_f>0$ rather than $M_f=1$.
As for the DNM above, our main goal is to estimate the complexity to enter the region
of quadratic convergence $\Q$ defined in \eqref{eq-quadr-conv}. Let us start
from some $x_0 \in \E$. Define the central path $x(t)$, $0
\leq t \leq 1$, by the following equation:
\beq\label{def-CP}
\ba{rcl}
\nabla f(x(t)) & = & t \nabla f(x_0).
\ea
\eeq
Clearly, $x(1) = x_0$ and $x(0) = x^*$. Note that this is a
trajectory of minimizers of the following parametric
family of general self-concordant functions:
\beq\label{eq-PFamily}
\ba{rcl}
x(t) & = & \arg\min\limits_{x \in \E} \left\{ \; f_t(x)
\Def f(x) - t \la \nabla f(x_0), x \ra \; \right\}, \quad
0 \leq t \leq 1.
\ea
\eeq
Let us introduce two constants
\beq\label{eq-Const}
\ba{rcl}
\beta & = & 0.026, \quad \gamma \; = \; 0.1125 < {\sqrt{\beta} \over
1 + \sqrt{\beta}} - \beta.
\ea
\eeq
We say that a point $x$ satisfies an {\em approximate
centering condition} if
\beq\label{eq-Approx}
\ba{rcl}
\lambda_{f_t}(x) & \equiv & \| \nabla f(x) - t \nabla f(x_0) \|^*_x \; \leq \; {\beta \over M_f}.
\ea
\eeq
Consider the path-following iterate: 
\beq\label{eq-ItPF}
(t_+,x_+) \; = \; {\cal P}(t,x) \; \Def \; \left\{
\ba{rcl}
t_+ & = & \max\left\{t - {\gamma \over M_f \| \nabla f(x_0) \|^*_x},0\right\},\\
\\
x_+ & = & x - [\nabla^2 f(x)]^{-1}(\nabla f(x) - t_+
\nabla f(x_0)).
\ea \right.
\eeq

The following statement is a counterpart of Theorem 4.2.8 in
\cite{nesterov2004introduction}.
\BT\label{th-ItPF}
If the pair $(x,t)$ satisfies (\ref{eq-Approx}), and $\beta$, $\gamma$ are chosen such that
\beq
\label{eq-ItPF-Appr-PF}
|\gamma| \leq \frac{\sqrt{\beta}}{1+\sqrt{\beta}} - \beta,
\eeq
then the pair $(x_+,t_+)$ satisfies (\ref{eq-Approx}) too.
\ET
\proof Let us denote $\lambda_0 = \| \nabla f(x) - t \nabla f(x_0) \|^*_x$, $\lambda_1 = \| \nabla f(x) - t_+ \nabla f(x_0) \|^*_x$ and $\lambda_+ = \| \nabla f(x_+) - t_+ \nabla f(x_0) \|^*_{x_+}$. Clearly, $\lambda_0 \leq \frac{\beta}{M_f}$. If $t_+>0$,  we have
\[
\lambda_1 = \left\| \nabla f(x) - t \nabla f(x_0) + \frac{\gamma}{M_f\| \nabla f(x_0) \|^*_x} \nabla f(x_0) \right\|^*_x \leq \frac{\beta+|\gamma|}{M_f}.
\] 
If $t_+=0$, we have that $t \leq {\gamma \over M_f \| \nabla f(x_0) \|^*_x}$. Hence,
\[
\lambda_1 =  \| \nabla f(x)\|^*_x=  \| \nabla f(x) - t \nabla f(x_0) + t \nabla f(x_0)\|^*_x \leq \frac{\beta+|\gamma|}{M_f}.
\]
Since $x_+$ is obtained from $x$ as the Standard Newton Step for the function $f_{t_+}(x)$, by \eqref{eq-StandQuad}, 
$
\lambda_+  \leq M_f \left( \frac{\lambda_1}{1-M_f\lambda_1} \right)^2. 
$ 
The statement of the theorem follows from the fact that inequality 
$ M_f\left( \frac{\lambda_1}{1-M_f\lambda_1} \right)^2 \leq \frac{\beta}{M_f}$ is equivalent to inequality $\lambda_1 \leq \frac{1}{M_f}\frac{\sqrt{\beta}}{1+\sqrt{\beta}}$. 
\qed

\begin{remark}
\label{RM:switch_to_Newton}
If at some iterate we obtain that $t_+=0$, then $t$ is not updated in the later iterates and the algorithm automatically switches to the standard full-step Newton method:
$
x_+   =   x - [\nabla^2 f(x)]^{-1}\nabla f(x).
$
According to \eqref{eq-StandQuad}, for this method it holds that 
\beq
\label{eq:switch_Newt_1}
\lambda_f(x_+)   \leq   \frac{M_f \lambda_f^2(x)}{(1-M_f\lambda_f(x))^2}.
\eeq
From Theorem \ref{th-ItPF}, we know that, if $t_+=0$, then the approximate centering condition \eqref{eq-Approx} holds at $(x_+,t_+)$, which means that, since $t_+=0$,
\beq
\label{eq:switch_Newt_2}
\lambda_f(x_+) = \| \nabla f(x_+)  \|^*_{x_+} \; \leq \; {\beta \over M_f}.
\eeq
Combining \eqref{eq:switch_Newt_1}, \eqref{eq:switch_Newt_2}, and our choice of the value $\beta$ in \eqref{eq-Const}, we see that the point $x_+$ belongs to the region of quadratic convergence of the standard Newton method. Thus, when $t_+=0$, the PFS automatically switches to the quadratically-convergent Newton method. 
\end{remark}

Let us prove the first main result of this section that gives convergence rate for the penalty parameter $t_k$ in the
PFS as applied to problem (\ref{prob-SCF}) when the region of quadratic convergence $\Q$ defined in \eqref{eq-quadr-conv} is not yet reached. 
The main novelty of this result is that we rely only on the assumption that $f$ is a self-concordant function rather than a self-concordant barrier as in \cite{nesterov2004introduction}.
\BT\label{th-PF}
Consider the path-following scheme (PFS):
\beq\label{eq-PF-SCF}
\ba{rcl}
t_0 = 1, \; x_0 \in \E, \quad (t_{k+1}, x_{k+1}) & = &
{\cal P}(t_k,x_k), \quad k \geq 0,
\ea
\eeq
where ${\cal P}$ is defined in \eqref{eq-ItPF}.
Assume that $\lambda_f(x_k) \geq {1 \over 2M_f}$ for all $k
= 0, \dots, N$. Then,
\beq\label{eq-PFRate}
\ba{rcl}
t_N & \leq & \exp \left\{ - {\gamma(\gamma -  2\beta)  N^2 \over 2M_f^2
(f(x_0) - f^*)} \right\}.
\ea
\eeq
\ET
\proof
Denote $c = - \nabla f(x_0)$. Then,
$$
\ba{rcl}
t_{k+1} & \refEQ{eq-ItPF} & t_k - {\gamma \over M_f \| c
\|^*_{x_k}} \; = \; t_k \left( 1 - {\gamma \over M_f t_k \|
c \|^*_{x_k}} \right) \; \leq \; t_k \exp \left\{ - {\gamma
\over M_f t_k \| c \|^*_{x_k}} \right\}.
\ea
$$
Thus, $t_N \leq \exp \left\{ - {\gamma \over M_f} S_N
\right\}$, where $S_N = \sum\limits_{k=0}^N {1 \over t_k
\| c \|^*_{x_k}}$.

Let us estimate the value $S_N$ from below. Note that
\beq\label{eq-Step}
\ba{rcl}
x_k - x_{k+1} & \refEQ{eq-ItPF} & [\nabla^2
f(x_k)]^{-1}\left(t_k c + \nabla f(x_k) - {\gamma c \over
M_f \| c \|^*_{x_k}} \right).
\ea
\eeq
Therefore,
\beq\label{eq-RK}
\ba{rcl}
r_k & \Def & \| x_k - x_{k+1} \|_{x_k} \; \refLE{eq-Approx} \; {\beta +
\gamma \over M_f}.
\ea
\eeq
On the other hand, ${\beta^2 \over M_f^2}
\refGE{eq-Approx} \lambda_f^2(x_k) + 2 t_k \la \nabla
f(x_k), [\nabla^2f(x_k)]^{-1} c \ra + t_k^2 (\| c
\|_{x_k}^*)^2$. Hence,
\beq\label{eq-Scal}
\ba{rcl}
- \la \nabla f(x_k), [\nabla^2f(x_k)]^{-1} c \ra & \geq &
{1 \over 2 t_k} \left[ \lambda_f^2(x_k) + t_k^2 (\| c
\|_{x_k}^*)^2 - {\beta^2 \over M_f^2}\right].
\ea
\eeq
Therefore, denoting $\lambda_k = \|\nabla f(x_k) - t_k \nabla f(x_0)\|_{x_k}^*$, 
\beq
\label{eq-th-PF-Proof-1}
\ba{l}
f(x_k) - f(x_{k+1}) \; \refGE{eq-UBound} \; \la \nabla
f(x_k), x_k - x_{k+1} \ra - {1 \over M_f^2} \omega_*(M_f
r_k)\\
\refEQ{eq-Step} \; \la \nabla f(x_k), [\nabla^2
f(x_k)]^{-1}\left(t_k c + \nabla f(x_k) - {\gamma c \over
M_f \| c \|^*_{x_k}} \right) \ra - {1 \over M_f^2}
\omega_*(M_f r_k)\\
= \lambda_k^2 - t_k \la c, [\nabla^2 f(x_k)]^{-1}(t_k c + \nabla f(x_k)) \ra + \la \nabla f(x_k), [\nabla^2 f(x_k)]^{-1} \left( \frac{-\gamma c}{M_f \| c \|^*_{x_k}} \right) \ra \\
\hspace{2em}- {1 \over M_f^2} \omega_*(M_f r_k)\\
\geq \lambda_k^2 - t_k\| c \|^*_{x_k} \lambda_k	- \frac{\gamma}{M_f \| c \|^*_{x_k}} \la \nabla f(x_k), [\nabla^2 f(x_k)]^{-1}c \ra - {1 \over M_f^2} \omega_*(M_f r_k)\\
\refGE{eq-Scal}  \lambda_k^2 - t_k\| c \|^*_{x_k} \lambda_k + \frac{\gamma}{2M_f t_k \| c \|^*_{x_k}} \left[ \lambda_f^2(x_k) + t_k^2 (\| c
\|_{x_k}^*)^2 - {\beta^2 \over M_f^2}\right] - {1 \over M_f^2} \omega_*(M_f r_k) \\
\refGE{eq-RK}  \frac{\gamma -2 M_f \lambda_k}{2M_f} t_k \| c \|^*_{x_k} + \rho_k \refGE{eq-Approx} \frac{\gamma -2\beta}{2M_f} t_k \| c \|^*_{x_k} + \rho_k,
\ea
\eeq
where $\rho_k = {\gamma \over 2 M_f t_k \| c \|^*_{x_k}}
\left[ \lambda_f^2(x_k)  - {\beta^2 \over M_f^2}\right] -
{1 \over M_f^2} \omega_*(\beta + \gamma)$.

Our next goal is to show that $\rho_k \geq 0$. Note that $t_k \| c \|^*_{x_k} \refLE{eq-Approx} \lambda_f(x_k) + \frac{\beta}{M_f}$.
Since $\lambda_f(x_k) \geq {1 \over 2M_f}$, we have
$$
\ba{rcl}
\rho_k & \geq & {\gamma \over 2 M_f }
\left[ \lambda_f(x_k)  - {\beta \over M_f}\right] -
{1 \over M_f^2} \omega_*(\beta + \gamma) \geq {\gamma (1-2 \beta) \over 4 M_f^2 } -
{1 \over M_f^2} \omega_*(\beta + \gamma).
\ea
$$
Using the values (\ref{eq-Const}), by direct computation
we can see that the right-hand side of this inequality is
positive.

Thus, we have proved that $f(x_k) - f(x_{k+1}) \geq {
\gamma -  2\beta  \over 2 M_f} t_k \| c \|^*_{x_k}$. Therefore,
\beq \label{eq-PF-SN-lower}
\ba{rcl}
S_N & \geq & \sum\limits_{k=0}^N {\gamma -  2\beta \over 2M_f
(f(x_k) - f(x_{k+1})) }\\
\\
& \geq & {\gamma -  2\beta \over 2M_f} \min\limits_{\tau \in
\R^{N+1}_+} \left\{ \sum\limits_{i=1}^{N+1} {1 \over
\tau^{(i)}}: \; \sum\limits_{i=1}^{N+1} \tau^{(i)} =
f(x_0) - f(x_{N+1}) \right\}\\
\\
& = & {(\gamma -  2\beta) (N+1)^2\over 2M_f(f(x_0)-f(x_{N+1}))} \geq {(\gamma -  2\beta) (N+1)^2\over 2M_f(f(x_0)-f^*)}.
\hfill{\mbox{\qed}}
\ea
\eeq

\begin{remark}
\label{RM:superlinear_rate}
The bound \eqref{eq-PFRate} can be slightly improved to give a superlinear convergence. Indeed, on the one hand, we have
\[
t_{k+1}  \refEQ{eq-ItPF}  t_k - {\gamma \over M_f \| c \|^*_{x_k}} \; = \; t_k \left( 1 - {\gamma \over M_f t_k \| c \|^*_{x_k}} \right).
\]
On the other hand, as we showed in the proof of Theorem \ref{th-PF}, $f(x_k) - f(x_{k+1}) \geq {\gamma -  2\beta  \over 2 M_f} t_k \| c \|^*_{x_k}$.
Whence, if $f(x_k) - f(x_{k+1}) \leq {\gamma(\gamma -  2\beta) \over 2M_f^2 }$, we have ${\gamma(\gamma -  2\beta) \over 2M_f^2 } \geq {\gamma -  2\beta  \over 2 M_f} t_k \| c \|^*_{x_k}$, and $t_k \leq {\gamma \over M_f\| c \|^*_{x_k}}$. Thus, recalling that $c=-\nabla f(x_0)$, from \eqref{eq-ItPF}, we have that $t_{k+1}=0$, and, according to Remark \ref{RM:switch_to_Newton}, the PFS automatically switches to the quadratically-convergent full-step Newton method.
Thus, since we are interested in the complexity of reaching the region of quadratic convergence, we assume that for $k=0,...,N$,  $f(x_k) - f(x_{k+1}) > {\gamma(\gamma -  2\beta) \over 2M_f^2 }$. Then, for $k=0,...,N$ we have
\[
t_{k+1} = t_k \left( 1 - {\gamma \over M_f t_k \| c \|^*_{x_k}} \right) \leq t_k \left( 1 - {\gamma(\gamma -  2\beta) \over 2M_f^2 (f(x_k) - f(x_{k+1}))} \right).
\]
Thus, denoting $\tau = {\gamma(\gamma -  2\beta) \over 2M_f^2 }$,  $\Delta_{k}=f(x_k) - f(x_{k+1})$, we obtain
\[
\ln t_N \leq \sum_{k=0}^{N}\ln\left(1-\frac{\tau}{\Delta_k}\right) \leq   \max_{\Delta_k \geq 0: \sum_{k=0}^{N} \Delta_k = f(x_0)-f(x_{N+1})} \sum_{k=0}^{N}\ln\left(1-\frac{\tau}{\Delta_k}\right) .
\]
The latter is a concave maximization problem. By symmetry, the minimum is achieved when $\Delta_k = \frac{ f(x_0)-f(x_{N+1})}{N+1} \leq \frac{ f(x_0)-f(x^*)}{N+1}$. So, we obtain
\beq
\label{eq:superlinear_rate}
t_N \leq  \left(1 - \frac{\gamma(\gamma -  2\beta) (N+1)}{2M_f^2(f(x_0)-f(x^*))}\right)^{N+1},
\eeq
which holds while $N+1 < \frac{f(x_0)-f(x^*)}{\tau}$.
Thus, we obtain that the PFS has \textit{global superlinear convergence}. We note also that the closer $f(x_0)$ to $f(x^*)$, the faster the algorithm converges. This is reasonable since if  $f(x_0)$ is close to $f(x^*)$, the point $x_0$ is close to the region of local quadratic convergence.
\end{remark}

Let us now estimate  the number of iterations, which is
sufficient for method (\ref{eq-PF-SCF}) to enter the region
of quadratic convergence $\Q$ defined in \eqref{eq-quadr-conv}. This is the second main result of this section. Denote 
$$
D = \max\limits_{x,y \in  {\rm dom} f} \{ \| x - y \|_{x_0} : \; f(x) \leq f(x_0), \; f(y) \leq f(x_0)  \}.
$$
\BT\label{th-CompPF}
Let sequence $\{ x_k \}_{k \geq 0}$ be generated by the
method (\ref{eq-PF-SCF}). Then, if
\beq\label{eq-CompPF}
\ba{rcl}
N & \geq &\left[\frac{2 \Delta(x_0)}{\gamma(\gamma - 2 \beta)} \ln
\frac{M_fD\omega^{-1}(\Delta(x_0))}{\omega\left(\frac{(1-\beta)(1-2\beta)}{2}  \right)} \right]^{1/2}
\ea
\eeq
we have $x_N \in \Q$.
\ET
\proof
Indeed,  
\[
f(x(t_k)) - f^* \leq \la \nabla f(x(t_k)), x(t_k) -
x^* \ra \refEQ{def-CP} t_k \la \nabla f(x_0), x(t_k) - x^*
\ra \leq t_k \lambda_f(x_0) D,
\]
where we used that $f(x_{k+1}) \leq f(x_k)$, $k \geq 0$, see \eqref{eq-th-PF-Proof-1}. 
Since $\omega(M_f
\lambda_f(x_0)) \refLE{eq-Decr} M_f^2 (f(x_0) - f^*)$, we have
$$
\ba{rcl}
{1 \over M_f^2}\omega(M_f \lambda_f(x(t_k))) \;
\refLE{eq-Decr} \; f(x(t_k)) - f^* & \leq & {t_k \over
M_f} \omega^{-1}(\Delta(x_0)) D.
\ea
$$
Note that $\| \nabla f(x_k) - \nabla f(x(t_k)) \|^*_{x_k}
\refEQ{def-CP} \| \nabla f(x_k) - t_k \nabla f(x_0)
\|^*_{x_k} \refLE{eq-Approx} {\beta \over M_f} < {1 \over M_f}$. Therefore,
$$
\ba{rcl}
\lambda_f(x_k) & \refLE{eq-Approx} & t_k \| \nabla f(x_0)
\|^*_{x_k} + {\beta \over M_f}   \refEQ{def-CP}   \la \nabla
f(x(t_k)), [\nabla^2
f(x_k)]^{-1} \nabla f(x(t_k)) \ra^{\frac{1}{2}} +{\beta \over M_f}\\
&\refLE{eq-CompD} & {1 \over 1 - \beta} \lambda_f(x(t_k))
+ {\beta \over M_f}.
\ea
$$
Thus, inclusion $x_k \in \Q$, is ensured by inequality
$\lambda_f(x(t_k)) \leq {(1-\beta)(1 - 2\beta) \over 2
M_f}$. Consequently, we need
\beq
\label{eq:t_k_for_quadratic_conv}
\ba{rcl}
{t_k \over M_f} \omega^{-1}(\Delta(x_0)) D & \leq & {1
\over M_f^2} \omega\left({(1-\beta)(1 - 2\beta) \over
2}\right) 
\ea
\eeq
It remains to use inequality (\ref{eq-PFRate}).
\qed


 As we can see from the estimate (\ref{eq-CompPF}) for the global complexity, up to a
logarithmic factor, the number of iterations of the
PFS is proportional to
$\Delta^{1/2}(x_0)$, where $\Delta(x_0)$ is defined in \eqref{eq-BoundCN}. This is much better than the
guarantee (\ref{eq-BoundCN}) for the DNM
(\ref{met-DNM}). However, as we will see in Section
\ref{sc-FStrong}, for some special subclasses of
self-concordant functions the performance estimate
(\ref{eq-CompPF}) can be significantly improved.

\begin{remark}
Note that, in the complexity bound \eqref{eq-CompPF}, the constant $\left[\frac{2}{\gamma(\gamma - 2 \beta)} \right]^{1/2} \leq 17.1$. The choice of the parameters $\beta$ and $\gamma$ is governed by the following aspects. First, from Theorem \ref{th-ItPF}, these parameters should satisfy \eqref{eq-ItPF-Appr-PF}. Second, $\rho_k$ in the proof of Theorem  \ref{th-PF} should be non-negative. Third, the complexity in \eqref{eq-CompPF} is proportional to $(\gamma(\gamma - 2 \beta))^{-1/2}$, which is desired to be as small as it is possible. This motivates the following maximization problem for optimal choice of $\beta$, $\gamma$.
\begin{align}
\max \gamma(\gamma - 2 \beta) \;\; \text{s.t.}& \label{eq-b-g-opt-1}  \\
\frac{\sqrt{\beta}}{1+\sqrt{\beta}} - \beta - \gamma &\geq 0 \label{eq-b-g-opt-2} \\
{\gamma (1-2 \beta) \over 4 } - \omega_*(\beta + \gamma) &\geq 0. \label{eq-b-g-opt-3} 
\end{align}
Figure \ref{fig-beta-gamma-plot} illustrates this optimization problem and the optimal objective value.

\begin{figure}[h!]
	\centering
	\includegraphics[width=0.5\linewidth]{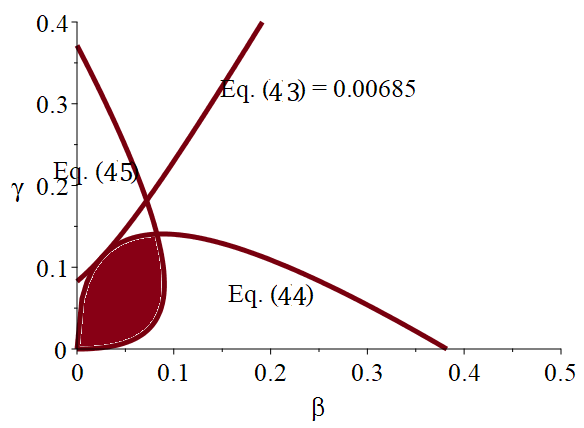}        
  \caption{Optimal choice of $\beta$ and $\gamma$. The feasible set given by \eqref{eq-b-g-opt-2} and \eqref{eq-b-g-opt-3} is filled with red color.}
	\label{fig-beta-gamma-plot}
\end{figure}
\end{remark}

%
%
%


The above results prescribe specific values for the accuracy of following the central path $\beta$ and the stepsize $\gamma$. Nevertheless, if we use a larger stepsize $\gamma$ and after the Newton step the approximate centering condition holds, we can continue to follow the path. This leads to a  PFS with adaptive choice of the stepsize $\gamma$ outlined below. 

\beq \notag
\ba{|c|}
\hline \\
\mbox{\bf Adaptive Path-Following Scheme}\\
\\
\hline \\
\quad \mbox{\MPT 
\begin{itemize}
	\item Set an initial point $x_0$, initial value of the penalty \\ parameter $t_0 = 1$, initial stepsize value $\gamma_{-1} \geq 0.1125$.
	\item $k$-th iteration. 
	Find the minimum value $i_k \geq 0$ s.t. the path-following step \eqref{eq-ItPF} with the stepsize $\gamma = 2^{1-i_k}\gamma_{k-1}$ outputs a point satisfying approximate centering condition \eqref{eq-Approx}.
	\item Set $\gamma_k = 2^{1-i_k}\gamma_{k-1}$, $x_{k+1} = x_+$, $t_{k+1} = t_+$.
\end{itemize}
\EPT}
\quad\\
\\
\hline
\ea
\eeq

By Theorem \ref{th-ItPF}, there exists $\hat{\gamma}$ s.t. the Newton step after the update of the parameter $t$ with the stepsize $\hat{\gamma}$ outputs a point satisfying the approximate centering condition. Hence, the search for $i_k$ is finite and $\gamma_k = 2^{1-i_k}\gamma_{k-1} \geq \frac{\hat{\gamma}}{2}$. Hence, the total number of Newton steps can be estimated as follows
\[
\sum_{j=0}^k i_k = \sum_{j=0}^k \left(1 + \log_2 \frac{\gamma_{j-1}}{\gamma_{j}} \right) = k+1 + \log_2 \frac{\gamma_{-1}}{\gamma_{k}} \leq k + \log_2 \frac{4\gamma_{-1}}{\hat{\gamma}}.
\] 
As we see, the price for the adaptivity is reasonable, taking into account that the practical performance of the adaptive algorithm is expected to be better since the penalty parameter $t$ potentially decreases faster.

Interestingly, increasing the stepsize $\gamma$ looks equivalent to decreasing the constant $M_f$. We underline that we assume the constant $M_f$ to be fixed and known. If the method would be adaptive also to the constant $M_f$, then the local estimate $M_f$ should be increased also in the inequality of the approximate centering condition. Also, note that for a fixed and known constant $M_f$ the stepsize for the DNM is optimal since it is obtained by the minimization of the upper bound for the self-concordant function. Thus, adaptive stepsize for this method does not make much sense for a fixed and known constant $M_f$.

%

\begin{remark}
After the first version of this paper appeared as the preprint \cite{dvurechensky2018global}, we were made aware of a concurrent work \cite{tran-dinh2018new}. In particular, the latter paper considers path-following methods for functions involving general self-concordant functions that are not barriers. At the same time there are several important differences between our work and theirs. They consider a different to our problem \eqref{prob-SCF} composite optimization problem of the form
\[
 \min\limits_{x \in \E} f(x)+g(x),
\]
where $f$ is a self-concordant function and $g$ is a simple closed convex function. Their algorithms and results are very specific to that problem and do not allow to obtain our algorithms and results as a particular case with $g=0$. First, their central path is defined by the inclusion
\[
0 \in t \nabla f(x(t)) - (1-t) g'(x_0) + \partial g(x(t)), 
\]
where $t \in [0,1]$, $\partial g(x)$ denotes the subdifferential of $g$,  $g'(x_0) \in \partial g(x_0)$. If $g=0$, then $g'(x_0)=0$ and the above equation, for any $t \in [0,1]$, is equivalent to the optimality condition for the problem \eqref{prob-SCF} and thus, following the path approximately does not make sense since for each $t \in [0,1]$ approximating $x(t)$ is equivalent to solving the original problem. In the same way, their algorithm is different to ours and when $g=0$, their main algorithmic step (12) is just the Standard Newton Step \eqref{eq-StandNStep} that may diverge and has only local convergence. Thus, for their results it is essential that $g \ne 0$ and our results do not follow from theirs as a particular case. Finally, they prove a linear convergence for their scheme, whereas we prove global superlinear convergence with explicit rate. We can say that their results and ours are complementary to each other. Finally, in the next sections we propose also a new predictor-corrector scheme, which they do not consider.
\end{remark}

z\section{Path-following scheme: implications for a feasibility problem}
\label{S:feasibility}

In this section, we first recall some properties of an important subclass of self-concordant functions, $\nu$-self-concordant barriers. These properties are then used to compare complexity of different methods applied to a feasibility problem which is solved by the minimization of a barrier on an affine subspace. Importantly, the proposed above PFS gives the best complexity, improving upon known complexities. Moreover, the obtained complexity corresponds to the \textit{global superlinear convergence} of the proposed scheme.

\subsection{Preliminaries on self-concordant barriers}
In this subsection, we recall the definition of self-concordant barriers and state their properties necessary to obtain the results of this section.

\BD
Let $F$ be a standard self-concordant function. We call it a $\nu$-self-concordant barrier for the set ${\rm Dom} F$ iff
\beq \label{eq-Def-SCB}
\sup_{y \in \E} [2 \la \nabla F(x),y\ra - \la \nabla^2 F(x) y,y\ra ] \leq \nu, 
\eeq
for all $x \in \dom F$. The value $\nu$ is called the parameter of the barrier.
\ED

Note that, if $\nabla^2 F(x)$ is non-degenerate, inequality \eqref{eq-Def-SCB} is equivalent to 
\beq\label{eq-Def-SCB-2}
\la [\nabla^2 F(x)]^{-1} \nabla F(x), \nabla F(x) \ra \leq \nu.
\eeq
The above inequality together with the duality relations \eqref{eq-Grad-Dual} and \eqref{eq-Hess-Dual} imply
\beq \label{eq-Dual-def-scb}
\la u , \nabla^2F_*(u)u \ra \leq  \nu , \quad u \in \dom F_*.
\eeq
\BT[Theorem 4.2.4 in \cite{nesterov2004introduction}]\label{th-scb-nf-dx-1}
1. Let $F(x)$ be a $\nu$-self-concordant barrier. Then, for any $x \in \dom F$ and $y \in {\rm Dom} F$, we have
\beq\label{eq-scb-nf-dx}
\la \nabla F(x), y-x \ra \leq \nu.
\eeq
2. A standard self-concordant function $F$ is a $\nu$-self-concordant barrier iff
\beq\label{eq-scb-log-charact}
F(y) \geq F(x) - \nu \ln \left( 1 - \frac{1}{\nu}\la \nabla F(x), y-x \ra \right).
\eeq
\ET


\subsection{Feasibility problem}
Let us consider the following feasibility problem
\beq
\label{eq:feasibility}
\text{Find} \;\; x \;\; \text{s.t.} \;\; x \in Q \;\; \text{and} \;\; Ax = b,
\eeq
where $x\in \R^n$, $b\in \R^m$, $A \in \R^{m\times n}$. We assume that $Q \subset \R^n$ is a closed and convex set and $0 \in \inter Q$. To solve this problem, we introduce a $\nu$-self-concordant \textit{barrier} $F(x)$ for the set $Q$ and minimize it over the affine manifold:
\beq
\label{eq:feasib_as_min}
\min_x F(x)  \;\; \text{s.t.} \;\; Ax=b.
\eeq
Without loss of generality, we assume that $0$ is the analytic center, i.e. $\nabla F(0) = 0$, and that  $ F(0) = 0$. Denote by $x^*$ a solution of \eqref{eq:feasib_as_min}. Then, clearly, $x^*$ solves the feasibility problem. We also assume that the \textit{feasibility depth}  $\e$ satisfies 
\beq
\label{eq:eps_def}
\e \Def \max \{\delta \geq 0 \; : \; (1-\delta)Q \cap \{x: Ax = b\}  \ne \emptyset \} >0.
\eeq
\BL\label{lm-feas_obj_resid}
Under the above assumptions, we have that
\begin{align}
& \frac{1}{\e}\geq 1 + \frac{1}{\nu}\la \nabla F(x^*),x^*\ra, \label{eq:feas_bound} \\
& F(x^*)-F(0) \leq \nu \ln\frac{1}{\e}. \label{eq:feas_obj_bound}
\end{align}
\EL
\proof
By the definition of the feasibility depth $\e$, there exists $\hat{x} \in (1-\e)Q$ such that $A\hat{x}=b$. Since $x^*$ is a solution to \eqref{eq:feasib_as_min}, there exists a Lagrange multiplier $y^*$ such that $ \nabla F(x^*) = A^T y^*$. Thus, for any $x$ satisfying  $Ax=b$, we have
\beq
\label{eq:y*b}
\la \nabla F(x^*),x\ra =  \la A^T y^*,x\ra =  \la y^*,Ax\ra = \la y^*,b\ra.
\eeq
Whence, $\la \nabla F(x^*),\hat{x}\ra = \la \nabla F(x^*),x^*\ra$. Thus, since $\hat{x} \in (1-\e)Q$, we have
\begin{align*}
\la \nabla F(x^*),x^*\ra &= \la \nabla F(x^*),\hat{x}\ra \leq \max_{u \in (1-\e)Q} \la \nabla F(x^*),u\ra \\
& = (1-\e) \max_{u\in Q} \{\la \nabla F(x^*),u-x^*\ra + \la \nabla F(x^*),x^*\ra \} \notag \\
&  \refLE{eq-scb-nf-dx} (1-\e)(\nu + \la \nabla F(x^*),x^*\ra).
\end{align*}
Thus, we have
\[
\e\la \nabla F(x^*),x^*\ra \leq (1-\e)\nu \Rightarrow \frac{1}{\nu}\la \nabla F(x^*),x^*\ra \leq \frac{1}{\e}-1,
\]
which is \eqref{eq:feas_bound}. Applying \eqref{eq-scb-log-charact} with $y=0$ and $x=x^*$, we obtain, by  \eqref{eq:feas_bound},
\[
F(0) \geq F(x^*) - \nu \ln \left( 1 + \frac{1}{\nu}\la \nabla F(x^*), x^* \ra \right) \geq F(x^*) - \nu \ln \frac{1}{\e},
\]
which is \eqref{eq:feas_obj_bound}.
\qed
Let us now consider the dual problem to \eqref{eq:feasib_as_min}. Introducing a Lagrange multiplier $y$ for the constraints $ A x = b$, we obtain
\[
\ba{rcl}
&\min_{x}&   \left\{ F(x) : \;\; A x = b \right\}  =  \min_{x} \max_y \{ F(x) + \la y, b - Ax\ra \} \\
\\
& = & \max_y  \la b, y\ra +  \min_x \{F(x) - \la A^Ty,x \ra  \} \} =
 \max_y \{\la b, y\ra - F_*(A^Ty) \}
\ea
\]
and the dual problem for \eqref{eq:feasib_as_min} is
\beq\label{eq-dual-feasib}
\min_{y} \left\{ \Phi(y) - \la b,y\ra \right\},
\eeq
where $\Phi(y) \Def F_*(A^Ty)$, $F_*$ being the Fenchel conjugate for $F$. We also have $\Phi(y^*) - \la b ,y^* \ra = -F(x^*)$, where $y^*$ is the solution of the dual problem \eqref{eq-dual-feasib}. Note that $\Phi(y) $ is a standard self-concordant function and is not a self-concordant barrier. Yet, it has a useful property \eqref{eq-Dual-def-scb} which allows to minimize it with the complexity standard for self-concordant barriers.

Assume now that we solve the dual problem \eqref{eq-dual-feasib} starting from $y=0$. Denote also $\widetilde{\Phi}(y) = \Phi(y) - \la b,y\ra$. Then, since $\widetilde{\Phi}(0)=\min_x F(x)=F(0)=0$, 
\beq
\label{eq:DeltaPhi_0}
\widetilde{\Phi}(0)-\widetilde{\Phi}^* = -\widetilde{\Phi}^* = F(x^*) = F(x^*) - F(0) \refLE{eq:feas_obj_bound} \nu \ln \frac{1}{\e}.
\eeq
We can apply the DNM \eqref{met-DNM} and the PFS \eqref{eq-PF-SCF} to solve the dual problem \eqref{eq-dual-feasib}.

Let us also consider dual path-following scheme \cite{nesterov1994interior}. To do so, consider 
the central path of a parametric family of problems
\beq
\label{eq:dualPFProbCehtral}
y_\sigma = \arg \min_y \{ \Phi(y): \sigma =  \la b, y \ra \}, \quad \sigma \geq 0.
\eeq
Since $\nabla \Phi(0) =A (\arg \min_x F(x))= 0$, we start with $y_0=0$ and $\sigma_0=0$, and the goal is to follow the path by increasing $\sigma$ until $\sigma=\sigma^*=\la b, y^*\ra$. Using \eqref{eq:y*b} and \eqref{eq:feas_bound}, we obtain
\beq
\label{eq:dualPFoptval}
\sigma^*=\la b, y^*\ra = \la \nabla F(x^*),x^*\ra \leq \nu\left(\frac{1}{\e} - 1\right) \leq \frac{\nu}{\e}.
\eeq
The Hessian of $\Phi$ induces  at $y$ local norm $\| \cdot\|_y$ and its conjugate $\| \cdot\|_y^*$. For the sake of simplicity, let us assume that we can follow the path \eqref{eq:dualPFProbCehtral} exactly. Namely, we use the following procedure
\beq\label{eq-ItDPF}
(\sigma_+,y_+) \; = \; {\cal DP}(\sigma,y) \; \Def \; \left\{
\ba{rcl}
\sigma_+ & = & \sigma + {\gamma \| b \|^*_y},\\
\\
y_+ & = & \arg \min_z \{ \Phi(z): \sigma_+ =  \la b, z \ra \}.
\ea \right.
\eeq
and iterate until $\sigma \geq \sigma^*$. We have
\[
\ba{rcl}
 \la b, y_+ \ra= \sigma_+  & = & \sigma+ \gamma \|b\|_y^* = \la b, y \ra + \gamma \|b\|_y^* = \la b, y \ra \left( 1+ {\gamma \|b\|_y^* \over \la b, y \ra }\right) \\
& \geq & \la b, y \ra \left( 1+ {\gamma \over \|y\|_y }\right) \geq \la b, y \ra \left( 1+ {\gamma \over \sqrt{\nu} }\right) = \sigma \left( 1+ {\gamma \over \sqrt{\nu} }\right) ,
\ea
\]
where we used that $\Phi$ is Fenchel conjugate for $\nu$-self-concordant barrier $F$  and, by \eqref{eq-Dual-def-scb}, $\|y\|_y = \la \nabla^2 \Phi(y)  y,  y \ra ^{1/2} \leq \sqrt{\nu}$ for all $y$. Hence, for the iterates $(\sigma_{k+1},y_{k+1}) = {\cal DP}(\sigma_k,y_k)$ it holds that $\sigma_k \geq \sigma_1  \left( 1+ {\gamma \over \sqrt{\nu} }\right)^k = \gamma \|b\|_{y_0}^*\left( 1+ {\gamma \over \sqrt{\nu} }\right)^k$ and, by \eqref{eq:dualPFoptval}, the complexity to get $\sigma_k \geq \sigma^*$ is bounded as
$
O \left(\frac{\sqrt{\nu}}{\gamma} \ln \frac{ \nu }{\gamma \e \|b\|_{y_0}^*} \right).
$

To sum up, we can apply three different strategies to solve problem \eqref{eq:feasib_as_min} by solving its dual \eqref{eq-dual-feasib}.
\begin{enumerate}
	\item {Damped Newton Method \eqref{met-DNM}.} In this case, we have that the complexity is given by \eqref{eq-BoundCN}, i.e. is $ O\left(\widetilde{\Phi}(0)-\widetilde{\Phi}^*\right) = O\left(\nu \ln \frac{1}{\e}\right)$.
	\item Path-following scheme \eqref{eq-PF-SCF}. In this case, we have that the complexity is given by \eqref{eq-CompPF}, i.e. is $O\left(\sqrt{\widetilde{\Phi}(0)-\widetilde{\Phi}^*}\right) = O\left(\sqrt{\nu \ln \frac{1}{\e}}\right)$.
	\item Dual path-following scheme \eqref{eq-ItDPF} has the complexity $O \left(\sqrt{\nu}\ln \frac{ \nu }{\e } \right)$. 
\end{enumerate}
Interestingly, the complexity $\sqrt{\nu \ln \frac{1}{\e}}$ corresponds to the superlinear convergence of the order $\exp\left(- \left(\frac{k}{\sqrt{\nu}}\right)^2\right)$, which is faster than the linear convergence $\exp\left(- \frac{k}{\sqrt{\nu}}\right)$ corresponding to the complexity $\sqrt{\nu} \ln \frac{1}{\e}$ that is typical for path-following methods. What is even more surprising is that we use the short-step PFS which was for a long time believed to be inferior to long-step methods having complexity $\sqrt{\nu} \ln \frac{1}{\e}$. Moreover to get the complexity $\sqrt{\nu \ln \frac{1}{\e}}$ we use the barrier property only to estimate the initial objective residual. We do not use the barrier property in the analysis of the method, only the self-concordance property is used.

Let us show that, actually, the feasibility problem \eqref{eq:feasibility} is very general. In particular, one can reduce a linear program in standard form to this feasibility problem. The primal-dual pair of linear programs corresponding to a linear program in standard form is
\beq
\label{eq:LP}
\min_x \{ \la c,x \ra: Ax =b, x \geq 0\} = \max_{s,y}\{ \la b,y \ra: s+A^Ty = c, s \geq0 \},
\eeq
where $x,s,c \in \R^n$, $y,b \in \R^m$, $A \in \R^{m\times n}$. If the matrix $A$ has full row rank $m$, without loss of generality, we can assume that the matrix $A$ has the form $(I_m,B)$, where $I_m$ is the identity matrix in $\R^{m \times m}$. Then, the vector $c$ can be divided in two blocks $(c_1,c_2)$ s.t. $c_1\in \R^m$, and, similarly, for the vector $s = (s_1,s_2)$.   Moreover, from the equation $s+A^Ty = c$, we have $s_1+y=c_1$, $s_2+B^Ty=c_2$. Hence, we have that the problem $\max_{s,y}\{ \la b,y \ra: s+A^Ty = c, s \geq0 \}$ is equivalent to $\max_{s}\{ \la b,c_1-s_1 \ra: s_2-B^Ts_1 = c_2-B^Tc_1, s \geq 0 \}$. By the strong duality, we have that solving this primal-dual pair of problems is equivalent to solving the system of equations and inequalities
\[
\la c,x \ra - \la b,c_1-s_1 \ra = 0, \; Ax =b, \; s_2-B^Ts_1 = c_2-B^Tc_1, \; x \geq 0, \; s \geq 0.
\]
Introducing a new variable $\tau \geq 0$, we have that the above system is equivalent to the following homogenous system of linear inequalities and equations
\[
\la c,x \ra + \la b, s_1 \ra - \tau\la b,c_1 \ra  = 0, \; Ax = \tau b, \; s_2-B^Ts_1 = \tau(c_2-B^Tc_1), \; x , s, \tau \geq 0.
\]
Thus, we obtain that solving the primal-dual pair of linear programs is equivalent to solving the system 
$
z \in \R^q, \; Qz = 0, \; z \geq 0,
$
where $z = (x,s,\tau)$ and the matrix $Q$ is composed of all the data in the previous system of equations (sometimes this representation is called {\em homogeneous self-dual embedding}, see \cite{ODonoghue2016conic} and references therein).
This system has a trivial solution $z=0$. To find a non-trivial solution, we add a constraint that $z$ belongs to the standard simplex $S_q(1) = \{z \in \R^q: z \geq 0, \la e_q, z\ra=1 \}$, where $e_q$ is the vector of all ones. This gives us a feasibility problem of finding $z$ s.t. $z \in S_q(1)$ and $Qz=0$. The only issue of this reformulation is that the standard simplex has empty interior. This can be easily resolved by transition to a $(q-1)$-dimensional problem. Indeed, by the simplex constraint, $z_q = 1-\sum_{i=1}^{q-1}z_i \geq 0$. Thus, if we introduce the variable $\bar{z}$ which consists of the first $q-1$ components of $z$, matrix $\bar{Q}$ which consists of the first $q-1$ columns of $Q$, vector $\bar{q}$ which is the last column of $Q$, we have that the considered feasibility problem is equivalent to 
$
\bar{z} \geq 0, \; \la e_{q-1},\bar{z}\ra \leq 1, \; (\bar{Q}-\bar{q}e_{q-1}^T)\bar{z} = -\bar{q}.
$

As we see, linear programming problem in the standard form \eqref{eq:LP} can be reformulated as the feasibility problem \eqref{eq:feasibility}. Thus, the obtained above complexity bounds for the feasibility problem apply also to linear programming problems.

\section{Minimizing self-concordant functions: predictor-corrector path-following scheme}\label{sc-PCPFS}

In this section, we return back to problem \eqref{prob-SCF} and propose an alternative algorithm, which we call predictor-corrector path-following scheme (PCPFS).
First, we prove an important technical lemma that the central path is sufficiently smooth. This allows us to use longer predictor steps in approximate tangent direction to the path, which is then followed by a corrector step to guarantee approximate centering condition. 
Notably, due to longer predictor step, we obtain better constant in the complexity bounds compared to the constant for the PFS \eqref{eq-PF-SCF}.

\subsection{Main lemma}


Given a fixed vector $c$, we consider the trajectory of minimizers of the following parametric
family of general self-concordant functions:
\beq\label{eq-PFamily1}
\ba{rcl}
x(t) & = & \arg\min\limits_{x \in \E} \left\{ \; f_t(x)
\Def f(x) + t \la c, x \ra \; \right\}, \quad
t \in \R.
\ea
\eeq
Recall that $\lambda_{f_t}(x)  \equiv \| \nabla f(x) + t c \|^*_x$, $t \in \R$.

The following lemma is the main result of this subsection and it will be used to justify our predictor-corrector path-following scheme (PCPFS).
\BL\label{lm-Main}
Let $f$ be an $M_f$-self-concordant function. Define $h = [\nabla^2 f(x)]^{-1} c $ and $r = \|h\|_x = \|c\|_x^*$. Then, for all $\tau$ such that $|\tau| M_f\|h\|_{x} < 1$,
\beq\label{eq-Main-Lemma}
\lambda_{f_{t-\tau}}(x+ \tau h ) \leq \lambda_{f_{t}}(x) \left(1+\frac{\tau M_f r }{1- \tau M_f r } \right) + \frac{1}{M_f}\left(\frac{ \tau M_f r }{1- \tau M_f r }\right)^2.
\eeq
\EL
\proof
Define $g_{\tau}=\nabla f(x+\tau h) + (t-\tau) c$ and $H_{\tau} = \nabla^2 f(x+\tau h)$. Then, for
\[
\ba{rcl}
\xi(\tau) & \Def & (\lambda_{f_{t-\tau}}(x+\tau h) )^2 = \la g_{\tau}, H_{\tau}^{-1} g_{\tau} \ra \\
					& = & \left\la  \nabla f(x+\tau h) + (t-\tau) c, [\nabla^2 f(x+\tau h)]^{-1} ( \nabla f(x+\tau h) + (t-\tau) c ) \right\ra,
\ea
\]
we have 
\beq\label{eq-Main-Lemma-proof1}
\ba{rcl}
\xi'(\tau) & = & - \left\la g_{\tau}, H_{\tau}^{-1} D^3f(x + \tau h)[h] H_{\tau}^{-1} g_{\tau} \right\ra + 2 \left\la H_{\tau}^{-1}(H_{\tau}h - c ) , g_{\tau} \right\ra \\
& \refLE{def-SCF-2} & 2M_f \|h\|_{x+\tau h} \|H_{\tau}^{-1} g_{\tau} \|_{x+\tau h}^2 + 2 \left\la  (H_{0}^{-1}- H_{\tau}^{-1})c , g_{\tau} \right\ra \\
& = & 2M_f \|h\|_{x+\tau h} \xi(\tau) + 2 \left\la  (H_{\tau}^{\frac{1}{2}}H_{0}^{-1}H_{\tau}^{\frac{1}{2}} - I)H_{\tau}^{-\frac{1}{2}}c, H_{\tau}^{-\frac{1}{2}}g_{\tau} \right\ra.
\ea
\eeq
Let us estimate $\|h\|_{x+\tau h}$ and the second term in the r.h.s.
Since $|\tau| M_f\|h\|_{x} < 1$, 
\beq
\label{eq-Main-Lemma-proof2}
\|h\|_{x+\tau h} = \frac{1}{\tau}\|(x+\tau h) - x \|_{x+\tau h} \refLE{eq-Norm-Compat} \frac{1}{\tau} \frac{\|(x+\tau h)-x\|_x}{1-M_f\|(x+\tau h)-x\|_x} = \frac{r}{1-\tau M_f r}.
\eeq
Further, since $|\tau| M_f\|h\|_{x} < 1$, from \eqref{eq-Compat}, we have
\[
\ba{rcl}
(1 - \tau M_f r)^2 H_0^{-1} & \preceq &
H_{\tau}^{-1} \; \preceq \; {1 \over (1 - \tau M_f r)^2} H_0^{-1}\\
\\
(1 - \tau M_f r)^2 H_{\tau}^{\frac{1}{2}} H_0^{-1} H_{\tau}^{\frac{1}{2}} & \preceq &
I \; \preceq \; {1 \over (1 - \tau M_f r)^2} H_{\tau}^{\frac{1}{2}}H_0^{-1}H_{\tau}^{\frac{1}{2}}\\
\\
(1 - \tau M_f r)^2 I & \preceq &
H_{\tau}^{\frac{1}{2}} H_0^{-1} H_{\tau}^{\frac{1}{2}} \; \preceq \; {1 \over (1 - \tau M_f r)^2} I \\
\\
\left((1 - \tau M_f r)^2 - 1\right) I  & \preceq &
H_{\tau}^{\frac{1}{2}} H_0^{-1} H_{\tau}^{\frac{1}{2}} - I \; \preceq \; \left({1 \over (1 - \tau M_f r)^2} - 1\right)  I. 
\ea
\]
Hence, 
\beq
\label{eq-Main-Lemma-proof3}
\ba{rcl}
&&\hspace{-2em}\left\la  (H_{\tau}^{\frac{1}{2}}H_{0}^{-1}H_{\tau}^{\frac{1}{2}}-I)H_{\tau}^{-\frac{1}{2}}c , H_{\tau}^{-\frac{1}{2}}g_{\tau} \right\ra  \leq  \| ( H_{\tau}^{\frac{1}{2}}H_{0}^{-1}H_{\tau}^{\frac{1}{2}}-I)H_{\tau}^{-\frac{1}{2}}c \| \|H_{\tau}^{-\frac{1}{2}} g_{\tau} \|\\
\\
&&  \leq  \left({1 \over (1 - \tau M_f r)^2} - 1 \right) \|H_{\tau}^{-\frac{1}{2}}c\| \sqrt{\xi(\tau)} \leq \left({1 \over (1 - \tau M_f r)^2} - 1 \right) {r \over 1 - \tau M_f r} \sqrt{\xi(\tau)},
\ea
\eeq
where we used that
\[
\ba{rcl}
\|H_{\tau}^{-\frac{1}{2}}c\|^2 & = & \left\la H_{\tau}^{-\frac{1}{2}}c, H_{\tau}^{-\frac{1}{2}}c \right\ra  \refLE{eq-Compat}  {1 \over (1 - \tau M_f r)^2} \left\la c, H_{0}^{-1} c \right\ra = {r^2 \over (1 - \tau M_f r)^2}.
\ea
\]
Combininig \eqref{eq-Main-Lemma-proof1}, \eqref{eq-Main-Lemma-proof2} and \eqref{eq-Main-Lemma-proof3}, we obtain
\[
\xi'(\tau) \leq {2M_fr \over 1 - \tau M_f r} \xi(\tau) + \left({1 \over (1 - \tau M_f r)^2} - 1 \right) {2r \over 1 - \tau M_f r} \sqrt{\xi(\tau)}.
\]
Using this inequality and denoting $\phi(\tau) = (1-\tau M_f r)^2 \xi(\tau)$, we obtain
\[
\ba{rcl}
\phi'(\tau) & = & -2M_f r(1-\tau M_f r) \xi(\tau) + (1-\tau M_f r)^2 \xi'(\tau)\\  
&  = & (1-\tau M_f r)^2 \left(\xi'(\tau) -  {2M_f r \over 1 - \tau M_f r}\xi(\tau)\right)\\
&  \leq & (1-\tau M_f r)^2 \left({1 \over (1 - \tau M_f r)^2} - 1 \right) {2r \over 1 - \tau M_f r} \sqrt{\xi(\tau)}\\
&  = &  2r \left({1 \over (1 - \tau M_f r)^2} - 1 \right) \sqrt{\phi(\tau) }.
\ea 
\]
Hence,
\[
\ba{rcl}
d \sqrt{\phi(\tau)} & \leq & \frac{1}{M_f}\left({1 \over (1 - \tau M_f r)^2} - 1 \right)  d (\tau M_f r).
\ea 
\]
Integrating both sides from $0$ to $\tau$, we obtain
\[
\ba{rcl}
\sqrt{\phi(\tau)} - \sqrt{\phi(0)}  & \leq & \frac{1}{M_f}\left({1 \over 1 - \tau M_f r} - \tau M_f r - 1 \right) = {\tau^2 M_f r^2 \over 1 - \tau M_f r}
\ea 
\]
\[
\text{and   } \qquad \qquad  \ba{rcl}
(1-\tau M_f r) \sqrt{\xi(\tau)} - \sqrt{\xi(0)}  & \leq &  {\tau^2 M_f r^2 \over 1 - \tau M_f r},
\ea 
\]
which finishes the proof since $\sqrt{\xi(\tau)}=\lambda_{f_{t-\tau}}(x+\tau h)$.
\qed

Let us make a remark on the meaning of the above Lemma. As we know from Section \ref{sc-PFS}, the quantity $\lambda_{f_{t}}(x)$ shows how good we follow the central path and this quantity should be kept sufficiently small during the iterations. The above lemma shows how this quantity evolves when we update the penalty parameter $t$. Unlike the PFS \eqref{eq-PF-SCF}, we will make simultaneously the update $t_+=t-\tau$ and $y_+=x+\tau h$ as the predictor step. The above lemma allows us to estimate the value $\lambda_{f_{t-\tau}}(x+ \tau h )$ and choose $\tau$ such that the corrector Newton step for $f_{t-\tau}$ again guarantees the approximate centering condition, but for the updated function $f_{t-\tau}$. This construction allows us to decrease $t$ faster leading to faster convergence.

\subsection{Predictor-corrector path-following scheme}
In this subsection, we describe our predictor-corrector path-following scheme (PCPFS) for solving problem (\ref{prob-SCF}). We again underline that we focus on the case when $f$ is a general self-concordant function without the barrier property. In the next section we will also show that the PCPFS has implications for problems with the barrier property of $f$.
As for the DNM and PFS above, our main goal in this section is to estimate the complexity to enter the region
of quadratic convergence $\Q$ defined in \eqref{eq-quadr-conv}.

For convenience, we recall from Section \ref{sc-PFS} the main objects important for our derivations. 
Namely, we start from some $x_0 \in \E$ and define the central path $x(t)$, $0
\leq t \leq 1$, by the following equation:
\beq\label{def-CP-PC}
\ba{rcl}
\nabla f(x(t)) & = & t \nabla f(x_0),
\ea
\eeq
or, equivalently,
\beq\label{eq-PFamily-PC}
\ba{rcl}
x(t) & = & \arg\min\limits_{x \in \E} \left\{ \; f_t(x)
\Def f(x) - t \la \nabla f(x_0), x \ra \; \right\}, \quad
0 \leq t \leq 1.
\ea
\eeq
In this section, we change the definition of the main parameters for following the central path (cf. \eqref{eq-Const} where the stepsize $\gamma$ is smaller)
\beq\label{eq-Const-PC}
\ba{rcl}
\beta & = & 0.0015, \quad \gamma \; = \; 0.158.
\ea
\eeq
As before, we say that point $x$ satisfies an {\em approximate
centering condition} if
\beq\label{eq-Approx-PC}
\ba{rcl}
\lambda_{f_t}(x) & \Def & \| \nabla f(x) - t \nabla
f(x_0) \|^*_x \; \leq \; {\beta \over M_f}.
\ea
\eeq
Instead of the PFS \eqref{eq-PF-SCF} we propose the predictor-corrector path-following iterate:
\beq\label{eq-ItPF-PC}
(t_+,x_+) \; = \; {\cal PC}(t,x) \; \Def \; \left\{
\ba{rcl}
t_+ & = & \max\left\{t - {\gamma \over M_f \| \nabla f(x_0) \|^*_x},0\right\},\\
\\
y		& = & x - {\gamma \over M_f \| \nabla f(x_0) \|^*_x}   [\nabla^2 f(x)]^{-1} \nabla f(x_0),\\
\\
x_+ & = & y - [\nabla^2 f(y)]^{-1}(\nabla f(y) - t_+ \nabla f(x_0)).
\ea \right.
\eeq
We refer to the $y$-step as the predictor and $x_+$-step as the corrector.

\BT\label{lm-ItPF-PC}
If the pair $(x,t)$ satisfies (\ref{eq-Approx-PC}), then the
pair $(x_+,t_+)$ satisfies (\ref{eq-Approx-PC}) too.
\ET
\proof
Applying Lemma \ref{lm-Main} with $c = -\nabla f(x_0)$, since, in this case, $r=\|h\|_x = \|\nabla f(x_0)\|_x^*$, we obtain 
\beq\label{eq-Lambda-f-t+-y}
\lambda_{f_{t_+}}(y) \leq \lambda_{f_{t}}(x) \left( 1+ \frac{\gamma}{1-\gamma}\right) + \frac{1}{M_f}\left(\frac{\gamma}{1- \gamma}\right)^2 \leq \frac{\beta}{M_f}\cdot\frac{1}{1-\gamma}  + \frac{1}{M_f}\left(\frac{\gamma}{1- \gamma}\right)^2.
\eeq
As $x_+$ is defined as the Standard Newton Step for $f_{t_+}$ from $y$, by \eqref{eq-StandQuad}, 
we have
\[
\lambda_{f_{t_+}}(x_+) \leq \frac{1}{M_f}\left(\frac{M_f\lambda_{f_{t_+}}(y)}{1-M_f\lambda_{f_{t_+}}(y)} \right)^2 = \frac{1}{M_f} \left(\omega'_*(M_f\lambda_{f_{t_+}}(y)) \right)^2.
\]
By the choice of $\beta$ and $\gamma$, we have $\omega'_*(M_f\lambda_{f_{t_+}}(y)) \leq \omega'_*\left(\frac{\beta}{1-\gamma}+ \left(\frac{\gamma}{1- \gamma}\right)^2 \right) \leq \sqrt{\beta}$, which finishes the proof.
\qed


Let us prove the first main result of this subsection that gives convergence rate for the penalty parameter $t_k$ in the
PCPFS as applied to problem (\ref{prob-SCF}) when the region of quadratic convergence $\Q$ defined in \eqref{eq-quadr-conv} is not yet reached. 
\BT\label{th-PF-PC}
Consider the predictor-corrector path-following scheme (PCPFS):
\beq\label{eq-PF-SCF-PC}
\ba{rcl}
t_0 = 1, \; x_0 \in \E, \quad (t_{k+1}, x_{k+1}) & = &
{\cal PC}(t_k,x_k), \quad k \geq 0,
\ea
\eeq
where ${\cal PC}$ is defined in \eqref{eq-ItPF-PC}.
Assume that $\lambda_f(x_k) \geq {1 \over 2M_f}$ for all $k
= 0, \dots, N$. Then
\beq\label{eq-PFRate-PC}
\ba{rcl}
t_N & \leq & \exp \left\{ - {\kappa \gamma N^2 \over M_f^2
(f(x_0) - f^*)} \right\}, \quad \kappa = \frac{\gamma}{2} - \frac{\beta}{(1-\gamma)^2} - \frac{\gamma^2}{(1-\gamma)^3}.
\ea
\eeq
\ET
\proof
Denote $c = - \nabla f(x_0)$. In the same way as in the proof of Theorem \ref{th-PF}, we obtain that
$t_N \leq \exp \left\{ - {\gamma \over M_f} S_N
\right\}$, where $S_N = \sum\limits_{k=0}^N {1 \over t_k
\| c \|^*_{x_k}}$.

Let us estimate the value $S_N$ from below. 
Note that 
\[
{\beta^2 \over M_f^2}
\refGE{eq-Approx-PC} \lambda_f^2(x_k) + 2 t_k \la \nabla
f(x_k), [\nabla^2f(x_k)]^{-1} c \ra + t_k^2 (\| c
\|_{x_k}^*)^2.
\]
Hence,
\beq\label{eq-Scal-PC}
\ba{rcl}
- \la \nabla f(x_k), [\nabla^2f(x_k)]^{-1} c \ra & \geq &
{1 \over 2 t_k} \left[ \lambda_f^2(x_k) + t_k^2 (\| c
\|_{x_k}^*)^2 - {\beta^2 \over M_f^2}\right].
\ea
\eeq
Thus, we obtain
\beq\label{eq-f-x-y}
\ba{l}
f(x_k) - f(y_k) \; \refGE{eq-UBound} \; 
\la \nabla f(x_k), x_k - y_k \ra - {1 \over M_f^2} \omega_*(M_f\|x_k-y_k\|_{x_k})\\
\refEQ{eq-ItPF-PC} \; - {\gamma \over M_f \|c \|^*_{x_k}} \la \nabla f(x_k), [\nabla^2 f(x_k)]^{-1}c \ra - {1 \over M_f^2}
\omega_*(\gamma)\\
\refGE{eq-Scal-PC}  \frac{\gamma}{2M_f t_k \| c \|^*_{x_k}} \left[ \lambda_f^2(x_k) + t_k^2 (\| c
\|_{x_k}^*)^2 - {\beta^2 \over M_f^2}\right] - {1 \over M_f^2}\omega_*(\gamma) \\
= \frac{\gamma}{2M_f} t_k \| c \|^*_{x_k} + \frac{\gamma}{2M_f t_k \| c \|^*_{x_k}} \left[ \lambda_f^2(x_k) - {\beta^2 \over M_f^2}\right] - {1 \over M_f^2} \omega_*(\gamma).
\ea
\eeq
Note that
\beq\label{eq-Step-Corr}
\ba{rcl}
y_k - x_{k+1} & \refEQ{eq-ItPF-PC} & [\nabla^2 f(y_k)]^{-1}\left(\nabla f(y_k) + t_{k+1}c\right).
\ea
\eeq
Hence,
\beq\label{eq-RK-PC}
\ba{rcl}
r_k & \Def & \| y_k - x_{k+1} \|_{y_k} \; = \lambda_{f_{t_{k+1}}}(y_k) \refLE{eq-Lambda-f-t+-y} \; \frac{1}{M_f}\left( \frac{\beta}{1-\gamma}+ \left(\frac{\gamma}{1- \gamma}\right)^2 \right).
\ea
\eeq
Therefore, 
\beq\label{eq-f-y-x}
\ba{l}
f(y_k) - f(x_{k+1}) \; \refGE{eq-UBound} \; \la \nabla
f(y_k), y_k - x_{k+1} \ra - {1 \over M_f^2} \omega_*(M_f
r_k)\\
\refEQ{eq-Step-Corr} \; \la \nabla f(y_k), [\nabla^2 f(y_k)]^{-1}\left(\nabla f(y_k) + t_{k+1} c \right) \ra - {1 \over M_f^2} \omega_*(M_f r_k)\\
=  \lambda_{f_{t_{k+1}}}^2(y_k) - \left(t_k - \frac{\gamma}{M_f \| c \|^*_{x_k}} \right)  \la c, [\nabla^2 f(y_k)]^{-1}\left(\nabla f(y_k) + t_{k+1} c  \right) \ra - {\omega_*(M_f r_k) \over M_f^2} \\
\geq \left(\lambda_{f_{t_{k+1}}}(y_k)\right)^2 - \left(t_k - \frac{\gamma}{M_f \| c \|^*_{x_k}} \right) \| c \|^*_{y_k} \lambda_{f_{t_{k+1}}}(y_k) - {1 \over M_f^2} \omega_*(M_f r_k)\\
\geq  -t_k\| c \|^*_{x_k} \frac{\lambda_{f_{t_{k+1}}}(y_k)}{1-\gamma} - {1 \over M_f^2} \omega_*(M_f r_k),
\ea
\eeq
where we used that $\|y_k-x_k\|_{x_k}\leq \gamma$ and, hence, $\| c \|^*_{y_k} \leq \frac{\| c \|^*_{x_k}}{1-\gamma}$ by \eqref{eq-Compat}.
Combining \eqref{eq-f-x-y} and \eqref{eq-f-y-x}, we obtain
\[
\ba{rcl}
f(x_k) - f(x_{k+1}) & \geq & t_k \| c \|^*_{x_k} \left(\frac{\gamma}{2M_f} - \frac{\lambda_{f_{t_{k+1}}}(y_k)}{1-\gamma}\right)
+  \frac{\gamma}{2M_ft_k \| c \|^*_{x_k}} \left[ \lambda_f^2(x_k) - {\beta^2 \over M_f^2}\right] \\
& & - {1 \over M_f^2} \omega_*(\gamma) - {1 \over M_f^2} \omega_*(M_f r_k) \\
&\refGE{eq-RK-PC}& \frac{\kappa}{M_f}  t_k \| c \|^*_{x_k} + \rho_k,
\ea
\]
where 
$
\kappa = \frac{\gamma}{2} - \frac{\beta}{(1-\gamma)^2} - \frac{\gamma^2}{(1-\gamma)^3},
$
$$
\rho_k = {\gamma \over 2 M_f t_k \| c \|^*_{x_k}}
\left[ \lambda_f^2(x_k)  - {\beta^2 \over M_f^2}\right] 
- {1 \over M_f^2} \omega_*(\gamma) - {1 \over M_f^2} \omega_*\left(\frac{\beta}{1-\gamma} - \frac{\gamma^2}{(1-\gamma)^2}\right).
$$

Our next goal is to show that $\rho_k \geq 0$. Note that $t_k \| c \|^*_{x_k} \refLE{eq-Approx-PC} \lambda_f(x_k) + \frac{\beta}{M_f}$.
Since $\lambda_f(x_k) \geq {1 \over 2M_f}$, we have
$$
\ba{rcl}
\rho_k & \geq & {\gamma \over 2 M_f }
\left[ \lambda_f(x_k)  - {\beta \over M_f}\right] -
{1 \over M_f^2} \omega_*(\gamma) - {1 \over M_f^2} \omega_*\left(\frac{\beta}{1-\gamma} - \frac{\gamma^2}{(1-\gamma)^2}\right) \\
& \geq & {\gamma (1-2 \beta) \over 4 M_f^2 } -
{1 \over M_f^2} \omega_*(\gamma) - {1 \over M_f^2} \omega_*\left(\frac{\beta}{1-\gamma} - \frac{\gamma^2}{(1-\gamma)^2}\right).
\ea
$$
Using the values (\ref{eq-Const-PC}), by direct computation
we can see that the right-hand side of this inequality is
positive.

Thus, we have proved that $f(x_k) - f(x_{k+1}) \geq \frac{\kappa}{M_f} t_k \| c \|^*_{x_k}$. Using the values \eqref{eq-Const-PC}, we can see that $\kappa >0$. 
Using the same steps as in the derivation of \eqref{eq-PF-SN-lower}, we obtain $S_N \geq {\kappa (N+1)^2\over M_f(f(x_0)-f(x_{N+1}))}\geq {\kappa (N+1)^2\over M_f(f(x_0)-f^*)}$.
\qed

Note that applying the same arguments as in Remark \ref{RM:superlinear_rate}, we see that the sequence $t_k$ for the PCPFS also has superlinear convergence.
Let us estimate now the number of iterations, which is
sufficient for PCPFS (\ref{eq-PF-SCF-PC}) to enter the region
of quadratic convergence $\Q$ defined in \eqref{eq-quadr-conv}. 
Recall that we denote $D = \max\limits_{x,y \in {\rm dom} f} \{ \| x - y \|_{x_0} : \; f(x) \leq f(x_0), \; f(y) \leq f(x_0)  \}$.
\BT\label{th-CompPF-PC}
Let sequence $\{ x_k \}_{k \geq 0}$ be generated by the
method (\ref{eq-PF-SCF-PC}). Then,  
\beq\label{eq-CompPF-PC}
\ba{rcl}
N & \geq &\left[{\Delta(x_0) \over \gamma \kappa }   \ln
\left({ M_fD \omega^{-1}(\Delta(x_0)) \over \omega\left({(1-\beta)(1 - 2\beta) \over 2} \right) } \right) \right]^{1/2}
\ea
\eeq
guarantees that $x_N \in \Q$.
\ET
\proof The result follows from (\ref{eq-PFRate-PC}) and \eqref{eq:t_k_for_quadratic_conv} that ensures that $x_k \in \Q$.
%
\qed

As we can see from the estimate (\ref{eq-CompPF-PC}), up to a
logarithmic factor, the number of iterations of the
PCPFS is proportional to
$\Delta^{1/2}(x_0)$. This is much better than the
guarantee (\ref{eq-BoundCN}) for the DNM \eqref{met-DNM}. 
Second, the constant $\left(\frac{1}{\gamma\kappa}\right)^{1/2} \leq  13.5$ in the complexity of the PCPFS is better than that of the standard PFS \eqref{eq-PF-SCF}, where it is about 17.

The above results prescribe specific values for the accuracy of following the central path $\beta$ and stepsize $\gamma$. Similarly to adaptive PFS in Section \ref{sc-PFS}, we can propose also an adaptive PCPFS.

\section{Predictor-corrector path-following scheme: implications for minimization problems}\label{sc-PCPFS-appl}
In this section, we apply the predictor-corrector path-following scheme (PCPFS) for minimization problems using an important subclass of self-concordant
functions, $\nu$-self-concordant barriers. The existing theory \cite{nesterov2004introduction} shows that path-following schemes have faster convergence when applied to self-concordant barriers. We show that our PCPFS as well possesses improved convergence properties if we use an additional barrier property. We first apply this approach to the classical minimization of a linear function over a set. Then, we consider applications to linearly constrained problems.

\subsection{Minimization problem: primal approach}
\label{S:min-primal}

In this subsection, we consider the following minimization problem
\beq\label{eq-Primal-Appr-Pr}
\min \, \la c, x \ra \quad \text{s.t.} \quad  x \in Q,
\eeq
where $Q = {\rm Dom}\, F$ is a bounded closed convex set with non-empty interior and $F$ is a $\nu$-self-concordant barrier meaning also that $M_F=1$. 
Similar problem was considered in \cite{nesterov2004introduction}, but we propose a different algorithm and obtain a better constant in the complexity bound. 
The reason is that we can make longer steps using the PCPFS.

To solve problem \eqref{eq-Primal-Appr-Pr}, we start from $x_0=x_F^*$, the analytic center of $Q$. 
For the purposes of this section we slightly redefine the main objects used in the path-following constructions.
We define the central path $x(t)$, $t \geq 0$, as
\beq\label{def-CP-Min}
\ba{rcl}
\nabla F(x(t)) & = & - t c.
\ea
\eeq
Clearly, $x(0) = x_0$. Note that 
\beq\label{eq-PFamily-Min}
\ba{rcl}
x(t) & = & \arg\min\limits_{x \in \E} \left\{ \; f_t(x)
\Def F(x) + t \la c, x \ra \; \right\}, \quad
t \geq 0.
\ea
\eeq
We again redefine the main parameters for following the central path (cf. \eqref{eq-Const} and \eqref{eq-Const-PC} where the stepsize $\gamma$ is smaller)
\beq\label{eq-Primal-Appr-Const}
\ba{rcl}
\beta & = & 0.06, \quad \gamma \; = \; 0.254.
\ea
\eeq
We say that point $x$ satisfies an {\em approximate
centering condition} if
\beq\label{eq-Approx-Min}
\ba{rcl}
\lambda_{f_t}(x) & \Def & \| \nabla F(x) + t c \|^*_x \; \leq \; {\beta }.
\ea
\eeq
Consider the following iterate which is a counterpart of \eqref{eq-ItPF-PC}:
\beq\label{eq-Primal-Appr-ItPF}
(t_+,x_+) \; = \; {\cal PC}(t,x) \; \Def \; \left\{
\ba{rcl}
t_+ & = & t + {\gamma \over \| c \|^*_x},\\
\\
y   & = & x - {\gamma \over \| c \|^*_x}   [\nabla^2 F(x)]^{-1} c,\\
\\
x_+ & = & y - [\nabla^2 F(y)]^{-1}(\nabla F(y) + t_+ c).
\ea \right.
\eeq

Similarly to Theorem \ref{lm-ItPF-PC}, we can prove that if the pair $(x,t)$ satisfies (\ref{eq-Approx-Min}), then the
pair $(x_+,t_+)$ satisfies (\ref{eq-Approx-Min}) too. Note that the chosen value $\gamma$ is nearly 2 times larger than the classical choice $\gamma = \frac{5}{36}$ in \cite{nesterov2004introduction}. This allows to make larger steps and the scheme converges faster.

\BL\label{lm-Primal-Appr-t-growth}
Let $x,t$ be generated by the iterates \eqref{eq-Primal-Appr-ItPF}. Then $t_+ \geq t \left( 1+ {\gamma \over \beta + \sqrt{\nu} }\right)$.
\EL
\proof
From \eqref{eq-Approx-Min} and \eqref{eq-Def-SCB-2}, we obtain
\[
t \|c\|_x^* = \|\nabla f_t(x) - \nabla F(x) \|_x^* \leq \|\nabla f_t(x) \|_x^*  + \|\nabla F(x) \|_x^* \leq \beta + \sqrt{\nu}
\]
Further, $t_+   =    t \left( 1 + {\gamma \over t \|c\|_x^* } \right) \geq t \left( 1 + {\gamma \over \beta + \sqrt{\nu} } \right)$.
\qed

\BL\label{lm-Primal-Appr-Approx}
Let $x^*$ be an optimal solution in \eqref{eq-Primal-Appr-Pr}. Then, for all $t > 0$,
\beq\label{eq-Primal-Appr-Obj-Central}
0 \leq \la c, x(t) \ra - \la c, x^* \ra \leq {\nu \over t  }.
\eeq
If $x$ satisfies approximate centering condition \eqref{eq-Approx-Min}, then
\beq\label{eq-Primal-Appr-Obj-Approx-Central}
|\la c, x \ra - \la c, x^* \ra | \leq   {1 \over t  } \left( \nu + {(\beta + \sqrt{\nu})\beta \over 1-\beta } \right).
\eeq
\EL
\proof
Since $x^*$ is an optimal solution in \eqref{eq-Primal-Appr-Pr}, by \eqref{def-CP-Min} and \eqref{eq-scb-nf-dx}
\[
0 \leq \la c, x(t) - x^* \ra = \frac{1}{t} \la \nabla F(x(t)), x^* -x(t) \ra \leq \frac{\nu}{t}.
\]
Further, let $x$ satisfy \eqref{eq-Approx-Min}. Denote $\lambda = \lambda_{f_t}(x)$. Then
\[
\ba{rcl}
t |\la c, x - x(t) \ra| & = & |\la \nabla f_t(x) - \nabla F(x), x - x(t) \ra | \\
& \leq & (\lambda + \sqrt{\nu}) \|x - x(t)\|_x \leq (\lambda + \sqrt{\nu})  \frac{\lambda}{1-\lambda} \leq {(\beta + \sqrt{\nu})\beta \over 1-\beta },
\ea
\]
where we used \eqref{eq-Def-SCB-2}, \eqref{eq-x-lambda-Bound}, and \eqref{eq-Approx-Min}. 
\qed

Note that inequality \eqref{eq-Primal-Appr-Obj-Approx-Central} gives us a reliable accuracy certificate based on the value of $t$.


The following is the main result of this subsection and gives complexity of the scheme \eqref{eq-Primal-Appr-ItPF} for solving problem \eqref{eq-Primal-Appr-Pr}.
\BT
Consider the predictor-corrector path-following scheme (PCPFS):
\beq\label{eq-PF-SCF-PC-Primal}
\ba{rcl}
t_0 = 0, \; x_0 \in {\rm dom} F \;\text{s.t.}\;  \|\nabla F(x_0)\|_{x_0}^*\leq \beta, \quad (t_{k+1}, x_{k+1}) & = &
{\cal PC}(t_k,x_k), \quad k \geq 0,
\ea
\eeq
where $\beta$ is defined in \eqref{eq-Primal-Appr-Const}, ${\cal PC}$ is defined in \eqref{eq-Primal-Appr-ItPF}.
Then, for any $\e >0$ and 
\[
N \geq N_{\e} = O \left( \frac{\sqrt{\nu}}{\gamma} \ln \left( \frac{\nu \|c\|_{x_F^*}^*}{\e }\right) \right),
\]
we have $x_N \in Q$ and $|\la c, x_N \ra- \la c, x^* \ra| \leq \e$.
\ET
\proof
By construction we have that $x_k \in Q$ for $k\geq 0$. 
Note that, by \eqref{eq-x-lambda-Bound}, 
$r_0 \Def \|x_0 - x_F^*\|_{x_0}  \leq \frac{\beta}{1-\beta}$. Hence, by \eqref{eq-Compat},
\[
\frac{\gamma}{t_1} = \|c\|_{x_0}^* \leq \frac{1}{1-r_0} \|c\|_{x_F^*}^* \leq \frac{1-\beta}{1-2 \beta} \|c\|_{x_F^*}^*.
\]
Hence, for all $k \geq 0$, from Lemma \ref{lm-Primal-Appr-t-growth}, $t_k \geq \frac{\gamma(1-2 \beta)}{(1-\beta)\|c\|_{x_F^*}^*} \left(1+\frac{\gamma}{\beta + \sqrt{\nu}} \right)^{k-1}$. 
As noted above, we have that \eqref{eq-Approx-Min} holds during the iterations. The result follows by the combination of \eqref{eq-Primal-Appr-Obj-Approx-Central} and the above lover bound for $t_k$.
\qed

As we can see, the main part of the obtained complexity bound is 
\[
3.94 \sqrt{\nu} \ln  \frac{\nu  \|c\|_{x_F^*}^*}{\e},
\]
which has a better constant 3.94 than the constant 7.2 for the one-step classical scheme \cite{nesterov2004introduction}. Similarly to the previous sections, the stepsize $\gamma$ may be chosen adaptively for better practical performance.

\subsection{Minimization problem: dual approach}
\label{S:min-dual}

In this subsection, we consider the following maximization problem
\beq\label{eq-Dual-Appr-Pr1}
\max - \la c, x \ra \quad \text{s.t.} \quad  B x = 0, \quad x \in Q,
\eeq
where $B \in \R^{m \times n}$ and $0 \in \inter Q$. The main difference with the previous subsection is that we now have linear constraints.
Let us introduce 
\[
A = \left(
\begin{aligned}
 - &c^T \\
& B\\
\end{aligned}
\right), \quad
b = \left(
\begin{aligned}
& 1 \\
& 0_m\\
\end{aligned}
\right).
\]
Then, problem \eqref{eq-Dual-Appr-Pr1} is equivalent to
\beq\label{eq-Dual-Appr-Pr2}
\max_{x,\alpha} \alpha \quad \text{s.t.} \quad A x = \alpha b, \quad x \in Q.
\eeq
Introducing penalty parameter $\sigma$ and a self-concordant barrier $F$ for the set $Q$ s.t. $\nabla F(0) = 0$, we consider the following parametric family of problems
\beq\label{eq-Dual-Appr-PF}
\max_{x,\alpha} \left\{ \alpha \sigma - F(x) : \quad A x = \alpha b \right\}, \quad \sigma \geq 0.
\eeq
Introducing a Lagrange multiplier $u$ for constraints $ A x = \alpha b$, we obtain
\[
\ba{rcl}
&& \max_{x,\alpha} \left\{ \alpha \sigma - F(x) : \quad A x = \alpha b \right\} =  \max_{x,\alpha} \min_u \{ \alpha \sigma - F(x) + \la u, Ax - \alpha b\ra \} \\
&& =  \min_u \{ \max_\alpha\{ \alpha (\sigma - \la u,b \ra) + \max_x \{-F(x) + \la A^Tu,x \ra  \} \} \\
&& =  \min_u \{F_*(A^Tu) : \sigma = \la b,u\ra \}
\ea
\]
and the dual problem for \eqref{eq-Dual-Appr-PF} is
\beq\label{eq-Dual-Appr-D1}
\min_{u} \left\{ \Phi(u) : \sigma = \la b,u\ra \right\}, \quad \sigma \geq 0,
\eeq
where $\Phi(u) \Def F_*(A^Tu)$, $F_*$ being the Fenchel conjugate for $F$. Note that $\Phi(u) $ is a standard self-concordant function with $M_{\Phi}=1$, but not a self-concordant barrier. 
Our goal is to follow the central path 
\[
u_\sigma \Def \arg \min_{u} \left\{ \Phi(u) : \sigma = \la b,u\ra \right\}
\] of the dual problem \eqref{eq-Dual-Appr-D1} for $\sigma$ from $0$ to $+\infty$. 
Optimality conditions in \eqref{eq-Dual-Appr-D1} give the following characterization of the dual central path
\beq\label{eq-Dual-Appr-PF-opt-cond}
\ba{rcl}
\nabla \Phi(u_\sigma) = A \nabla F_*(A^Tu_\sigma) & = & \alpha_\sigma b, \\
\la b , u_\sigma \ra & = & \sigma,
\ea
\eeq
where $u_\sigma$ is an optimal solution for \eqref{eq-Dual-Appr-D1} and $\alpha_\sigma$ is an optimal Lagrange multiplier. Let us define $x_\sigma \Def \nabla F_*(A^Tu_\sigma)$. Since $F_*$ is the Fenchel conjugate for a barrier $F$ for $Q$, $x_\sigma \in Q$ and $(\alpha_\sigma,x_\sigma)$ is a feasible point for \eqref{eq-Dual-Appr-Pr2}, see \eqref{eq-Grad-Dual}.

We use the same $\beta,\gamma$ as in the previous subsection (cf. \eqref{eq-Primal-Appr-Const}):
\beq\label{eq-Dual-Appr-Const}
\ba{rcl}
\beta & = & 0.06, \quad \gamma \; = \; 0.254.
\ea
\eeq
For the dual problem \eqref{eq-Dual-Appr-D1}, we define $\Phi_t (u) \Def   \Phi(u) - t b$ and say that a point $u$ satisfies an {\em approximate
centering condition} if, for some $t \in \R$,
\beq\label{eq-Dual-Appr-Approx-Min}
\ba{rcl}
\lambda_{\Phi_t}(u) & \Def & \| \nabla \Phi(u) - t b \|^*_u \; \leq \; {\beta },
\ea
\eeq
where the local norm $\|\cdot\|_u$ is induced by the Hessian of $\Phi(u)$.
We also define $t(u) \Def \arg \min_t \lambda_{\Phi_t}(u) $ and $\lambda(u) \Def \lambda_{\Phi_{t(u)}}(u) $.
Note that, for any $v$
\beq
\label{eq-t-of-v}
\ba{rcl}
t(v) &=& \arg \min_t \lambda_{\Phi_t}(v) = \arg \min_t (\| \nabla \Phi(v) - t b \|^*_v)^2 \\
&=& \arg \min_t \{\lambda_{\Phi}(v)^2 -2t \la b, [\nabla^2 \Phi(v)]^{-1} \nabla \Phi(v) \ra + t^2(\|b\|_v^*)^2 \}\\
&=& {\la b, [\nabla^2 \Phi(v)]^{-1} \nabla \Phi(v) \ra \over (\|b\|_v^*)^2 }.
\ea
\eeq

We consider the dual PCPFS that approximately follows the central path of the dual problem \eqref{eq-Dual-Appr-D1} in order to solve the primal problem \eqref{eq-Dual-Appr-Pr1}
\beq\label{eq-Dual-Appr-ItPF}
(\sigma_+,u_+) \; = \; {\cal PC}(\sigma,u) \; \Def \; \left\{
\ba{rcl}
v & = & u + {\gamma \over \| b \|^*_u}   [\nabla^2 \Phi(u)]^{-1} b, \\
\\
\sigma_+ & = & \la b, v \ra,\\
\\
u_+ &=& v+ \arg \min_{h: \la b, h \ra =0  }\left\{ \la \nabla \Phi(v), h\ra + \frac{1}{2} \la \nabla^2 \Phi(v) h,h\ra  \right\}.\\
\ea \right.
\eeq
Note that by the optimality condition for the minimization problem defining $u_+$, there exist a Lagrange multiplier $\kappa$ s.t.
\[
\nabla \Phi(v) + \nabla^2 \Phi(v)h  - \kappa b = 0 \quad \Leftrightarrow \quad h= - [\nabla^2 \Phi(v)]^{-1} (\nabla \Phi(v) - \kappa b).
 \]
Substituting this into the equality constraint $ \la b, h \ra =0 $, we obtain
\beq
\label{eq-kappa}
\kappa =  {\la b, [\nabla^2 \Phi(v)]^{-1} \nabla \Phi(v) \ra \over (\|b\|_v^*)^2 } \refEQ{eq-t-of-v} t(v).
\eeq
Moreover, 
\beq
\label{eq-u-step}
u_+ = v - [\nabla^2 \Phi(v)]^{-1} (\nabla \Phi(v) - t(v) b).
\eeq

\BL\label{lm-Dual-Appr-ItPF}
If $\lambda(u) \leq \beta$, then $\lambda(u_+) \leq \beta$ too.
\EL
\proof
Note that Lemma \ref{lm-Main} holds for all $t$. Thus, applying Lemma \ref{lm-Main} for $t=t(u)$ with $c = b$, $h = [\nabla^2 \Phi(u)]^{-1} b$, $\tau = {\gamma \over \|b\|_u^*}$, and $\Phi$ with $M_{\Phi}=1$ instead of $f$, since, in this case, $r=\|h\|_u = \|b\|_u^*$, we obtain, by \eqref{eq-Dual-Appr-ItPF} that $v=u+\tau h$ and  
\[
\lambda (v) = \lambda_{\Phi_{t(v)}}(v) = \min_t \lambda_{\Phi_{t}}(v) \leq \lambda_{\Phi_{t(u)-\tau}}(v) \leq  \frac{\lambda(u)}{1-\gamma} + \left(\frac{\gamma}{1- \gamma}\right)^2 \leq  \frac{\beta}{1-\gamma} + \left(\frac{\gamma}{1- \gamma}\right)^2. 
\]
Since, by \eqref{eq-u-step}, $u_+$ is obtained by the Standard Newton Step for $\Phi_{t(v)}$ from $v$, by \eqref{eq-StandQuad}, we have
\[
\lambda (u_+)\leq \lambda_{\Phi_{t(v)}}(u_+) \leq \left(\frac{\lambda_{\Phi_{t(v)}}(v)}{1-\lambda_{\Phi_{t(v)}}(v)} \right)^2 =  \left(\omega'_*(\lambda_{\Phi_{t(v)}}(v)) \right)^2.
\]
By the choice of $\beta$ and $\gamma$ in \eqref{eq-Dual-Appr-Const} and the estimate for $\lambda_{\Phi_{t(v)}}(v)$ above, we have $\omega'_*(\lambda_{\Phi_{t(v)}}(v)) \leq \omega'_*\left(\frac{\beta}{1-\gamma}+ \left(\frac{\gamma}{1- \gamma}\right)^2 \right) \leq \sqrt{\beta}$, which finishes the proof.
\qed

\BL\label{lm-Dual-Appr-sigma-growth}
Let $u,\sigma$ be generated by the iterates \eqref{eq-Dual-Appr-ItPF}. Then $\sigma_+ \geq \sigma \left( 1+ {\gamma \over \sqrt{\nu} }\right)$.
\EL
\proof
We have
\[
\ba{rcl}
\sigma_+ =  \la b, u_+ \ra & = & \la b, v \ra = \la b, u \ra + \gamma \|b\|_u^* = \la b, u \ra \left( 1+ {\gamma \|b\|_u^* \over \la b, u \ra }\right) \\
& \geq & \la b, u \ra \left( 1+ {\gamma \over \|u\|_u }\right) \geq \la b, u \ra \left( 1+ {\gamma \over \sqrt{\nu} }\right) = \sigma \left( 1+ {\gamma \over \sqrt{\nu} }\right) ,
\ea
\]
where we used that $\Phi$ is Fenchel conjugate for $\nu$-self-concordant barrier $F$  and, by \eqref{eq-Dual-def-scb}, $\|u\|_u = \la \nabla^2 \Phi(u)  u,  u \ra ^{1/2} \leq \sqrt{\nu}$ for all $u$.
\qed

\BL\label{lm-Dual-Appr-Primal-Approx}
Let $\alpha^*$ be an optimal function value in \eqref{eq-Dual-Appr-Pr2}. Then, for all $\sigma \geq 0$,
\beq\label{eq-Dual-Appr-Primal-Central}
0 \leq \alpha^* - \alpha_\sigma \leq {\nu \over \sigma  }.
\eeq
 If $u$ satisfies $\lambda(u) \leq \beta$ and $\la b, u \ra = \sigma$, then
\beq\label{eq-Dual-Appr-Primal-Approx-Central}
|\alpha^* - t(u)| \leq   {1 \over \sigma  } \left( \nu + {2\beta(1 - \beta) \over 1-2\beta } \sqrt{\nu} \right).
\eeq
\EL
\proof
Since $(\alpha^*,x^*)$ is an optimal solution for \eqref{eq-Dual-Appr-Pr2}, $A(x^*-x_\sigma) = (\alpha^*-\alpha_\sigma)b$ and
\begin{align*}
0 &\leq \alpha^* - \alpha_\sigma  =  {\la A(x^*-x_\sigma),u_\sigma \ra \over \la b , u_\sigma \ra } = {\la x^*-x_\sigma, A^Tu_\sigma \ra \over \la b , u_\sigma \ra } = {\la x^*-x_\sigma, \nabla F(x_\sigma) \ra \over \la b , u_\sigma \ra } \\
& \stackrel{\eqref{eq-scb-nf-dx}}{\leq} {\nu \over \la b , u_\sigma \ra } = {\nu \over \sigma  },
\end{align*}
where we used that $x_\sigma = \nabla F_*(A^Tu_\sigma)$ and, hence, $A^T u_\sigma  = \nabla F(x_\sigma)$.

Since $u_\sigma$ is a minimizer of self-concordant function $\Phi(v)-t(u)\la b,v\ra$ on the affine set $\la b,v \ra = \sigma$, by \eqref{eq-x-lambda-Bound}, we have
\[
\|u_\sigma-u\|_u \leq \omega_*'\left(\lambda_{\Phi(\cdot)-t(u)\la b,\cdot\ra}(u)\right)=\omega_*'\left(\lambda(u)\right) \leq \frac{\beta}{1-\beta}.
\]
Then, using Lemma \ref{lm-scf-grad-diff}, we obtain
\[
\ba{rcl}
&&\hspace{-2em}|\alpha_\sigma - t(u)| \sigma  =   |\alpha_\sigma - t(u)| \la b , u \ra   =  {|\la \nabla \Phi(u_\sigma)-\nabla \Phi(u) + \nabla \Phi(u) - t(u) b,u \ra| } \\
& \leq & |\la \nabla \Phi(u_\sigma)-\nabla \Phi(u),u \ra| + {|\la\nabla \Phi(u) - t(u) b,u \ra|  }\\
& \leq & { \|\nabla \Phi(u_\sigma)-\nabla \Phi(u)\|_u^* \|u\|_u } + { \|\nabla \Phi(u) - t(u) b\|_u^* \|u\|_u } \\
&\leq &
{ \|u_\sigma-u\|_u \over 1- \|u_\sigma-u\|_u } {\|u\|_u } + { \|\nabla \Phi(u) - t(u) b\|_u^* \|u\|_u } 
\leq  \left({\beta \over 1-2\beta} + \beta \right) \sqrt{\nu} = {2\beta(1 - \beta) \over 1-2\beta } \sqrt{\nu},
\ea
\]
where we used that $\Phi$ is Fenchel conjugate for $\nu$-self-concordant barrier $F$  and, by \eqref{eq-Dual-def-scb}, $\|u\|_u = \la \nabla^2 \Phi(u)  u,  u \ra ^{1/2} \leq \sqrt{\nu}$ for all $u$.
\qed
Note that inequality \eqref{eq-Dual-Appr-Primal-Approx-Central} gives us a reliable accuracy certificate based on the value of $\sigma$.
%

The following is the main result of this subsection and gives complexity of the scheme \eqref{eq-Dual-Appr-ItPF} for solving problem \eqref{eq-Dual-Appr-Pr1} in its equivalent form \eqref{eq-Dual-Appr-Pr2}.
\BT
Consider the dual predictor-corrector path-following scheme:
\beq\label{eq-PF-SCF-PC-Dual}
\ba{rcl}
\sigma_0=0, u_0=0, \quad (\sigma_{k+1}, u_{k+1}) & = &
{\cal PC}(\sigma_k,u_k), \quad k \geq 0,
\ea
\eeq
where ${\cal PC}$ is defined in \eqref{eq-Dual-Appr-ItPF}. Let also, for $k\geq 0$, $x_{k} = \nabla F_*(A^Tu_{k})$ and  $\hat{x}_k$ be the solution to the following problem
	\[
		\min \left\{\frac{1}{2}\|x-x_k\|_{x_k}^2 \;\; \text{ s.t. } \;\; Ax=t(u_k)b  \right\}.
	\]

Then, $\hat{x}_{k} \in Q$, $A\hat{x}_{k} = t(u_{k})b$ for $k\geq 0$. Moreover, for any $\e >0$ and 
\beq\label{eq-PF-SCF-PC-Dual-N}
N \geq N_{\e} = O \left( \frac{\sqrt{\nu}}{\gamma} \ln \left( \frac{\nu}{\e \gamma \|b\|_0^*}  \left( 1+ {2\beta(1 - \beta) \over (1-2\beta) \sqrt{\nu} }  \right)\right) \right),
\eeq
we have $|\alpha^* - t(u_{N})| \leq \e$.
\ET
\proof
By construction, we have $A\hat{x}_{k} = t(u_{k})b$ and our next step is to show that $\hat{x}_{k} \in Q$. We do that by showing that $\|\hat{x}_{k}-x_k\|_{x_k} \leq \beta <1$, which since $x_k \in Q$ implies that $\hat{x}_{k}$ is in the Dikin ellipsoid around $x_k$ and  belongs to $Q$. 

By the optimality conditions for the problem defining $\hat{x}_{k}$, we have that there exist a dual variable $y$ such that
\[
\ba{rcl}
\nabla^2 F(x_k)(\hat{x}_{k}-x_k)+A^Ty & = & 0, \\
A \hat{x}_{k} & = & t(u_k) b.
\ea
\]
From the first equality, we obtain 
$\hat{x}_{k}-x_k = - [\nabla^2 F(x_k)]^{-1} A^T y$.
Then, from the second equality, we obtain
$
y = \left(A [\nabla^2 F(x_k)]^{-1} A^T\right)^{-1} (A x_k-t(u_k) b),
$
which finally gives
\beq
\label{eq:DualPFmainTh_proof_2}
\hat{x}_{k}-x_k = - [\nabla^2 F(x_k)]^{-1} A^T \left(A [\nabla^2 F(x_k)]^{-1} A^T\right)^{-1} (A x_k-t(u_k) b).
\eeq
Recall that $A x_k= A \nabla F_*(A^Tu_{k}) = \nabla \Phi(u_k)$. Further, since $\Phi(u) = F_*(A^Tu)$, by \eqref{eq-Hess-Dual},
we have that $\nabla^2 \Phi(u_k) = A [\nabla^2 F(x_k)]^{-1} A^T$.   
Thus, from  \eqref{eq:DualPFmainTh_proof_2}, we have that
$
\hat{x}_{k}-x_k = - [\nabla^2 F(x_k)]^{-1} A^T [\nabla^2 \Phi(u_k)]^{-1}(\nabla \Phi (u_k)-t(u_k) b).
$
Hence, 
\[
\ba{rcl}
&&\hspace{-1em}\|\hat{x}_{k}-x_k\|_{x_k}^2  =  \la [\nabla^2 F(x_k)] [\nabla^2 F(x_k)]^{-1} A^T [\nabla^2 \Phi(u_k)]^{-1}(\nabla \Phi (u_k)-t(u_k) b), \\ 
& & [\nabla^2 F(x_k)]^{-1} A^T [\nabla^2 \Phi(u_k)]^{-1}(\nabla \Phi (u_k)-t(u_k) b) \ra \\
& = & \la  A [\nabla^2 F(x_k)]^{-1} A^T  [\nabla^2 \Phi(u_k)]^{-1}(\nabla \Phi (u_k)-t(u_k) b), [\nabla^2 \Phi(u_k)]^{-1}(\nabla \Phi (u_k)-t(u_k) b) \ra \\
& = & \la    \nabla \Phi (u_k)-t(u_k) b, [\nabla^2 \Phi(u_k)]^{-1}(\nabla \Phi (u_k)-t(u_k) b) \ra \\
& = & (\|\nabla \Phi (u_k)-t(u_k) b\|_{u_k}^*)^2 \leq \beta^2,
\ea
\]
where in the last inequality we used the approximate centering condition $\lambda(u_k)\leq \beta$. 
Thus, we have that $\hat{x}_{k}$ belongs to the Dikin ellipsoid around $x_k$ and, hence, belongs to $Q$.

Finally, we show the complexity result. 
Since $\nabla F(0) = 0$, we have that $\nabla \Phi(0) = 0$ and $\lambda(u_0) = \lambda(0) = 0 \leq \beta$. Hence, by Lemma \ref{lm-Dual-Appr-ItPF}, we have that $\lambda(u_k)\leq \beta$ for all $k\geq 0$. By \eqref{eq-Dual-Appr-ItPF} we have that $\sigma_k=\la b, u_k \ra$. Thus, by Lemma \ref{lm-Dual-Appr-Primal-Approx}, we have that \eqref{eq-Dual-Appr-Primal-Approx-Central} holds with $\sigma=\sigma_k$, $u=u_k$ for all $k\geq 0$.
Further, $\sigma_1 = \gamma \|b\|_0^*$ and, by Lemma \ref{lm-Dual-Appr-sigma-growth}, $\sigma_k \geq \gamma \|b\|_0^* \left(1+{\gamma \over \sqrt{\nu} }\right)^{k-1}$. 
Combining these observations with \eqref{eq-PF-SCF-PC-Dual-N}, we obtain that $|\alpha^* - t(u_{N})| \leq \e$.
\qed

As we can see, the main part of the obtained complexity bound is 
\[
3.94 \sqrt{\nu} \ln  \frac{\nu}{\e \|b\|_0^*},
\]
which has a better constant 3.94 than the constant 7.2 for the one-step classical scheme \cite{nesterov2004introduction}. Similarly to the previous sections, the stepsize $\gamma$ may be chosen adaptively for better practical performance. Note also that it is not necessary to calculate the point $\hat{x}_{k}$ in each iteration.

\section{Minimizing strongly convex functions with Lipschitz Hessian}
\label{sc-FStrong}

Let $B = B^* \succ 0$ map $\E$ to $\E^*$. Define
the Euclidean metric $\| x \|^2   =   \la B x, x \ra^{1/2},   x \in \E$.
In this section, we return back to the unconstrained problem \eqref{prob-SCF}, 
where $f$ is a strongly convex function with parameter $\sigma_f > 0$:
\beq\label{eq-strong}
\ba{rcl}
f(y ) & \geq & f(x) + \la \nabla f(x), y - x \ra + \half
\sigma_f \| y - x \|^2, \quad x, y \in \E.
\ea
\eeq
We also assume that  $f \in \C^3(\E)$ and its Hessian is Lipschitz
continuous:
\beq\label{eq-Hess}
\ba{rcl}
\| \nabla^2 f(x) - \nabla^2 f(y) \| & \leq H_f \| x - y
\|, \quad x, y \in \E.
\ea
\eeq
\BL\label{lb-SS}
Under assumptions above, function $f$ is self-concordant
with
\beq\label{eq-SS}
\ba{rcl}
M_f & = & {H_f \over 2 \sigma_f^{3/2}}.
\ea
\eeq
\EL
\proof
Indeed, for any point $x \in \E$ and direction $h \in \E$
we have
$$
\ba{rcl}
D^3f(x)[h]^3 & \refLE{eq-Hess} & H_f \left[ \| h \|^2
\right]^{3/2} \; \refLE{eq-strong} \; H_f \left[ {1 \over
\sigma_f} \la \nabla^2 f(x) h, h \ra \right]^{3/2}.
\ea
$$
It remains to use definition (\ref{def-SCF}).
\qed

Thus, problem (\ref{prob-SCF}) can be solved by methods
(\ref{met-DNM}) and (\ref{eq-PF-SCF}). The corresponding
complexity bounds can be given in terms of the complexity
measure $\Delta(x_0)   =   {H_f^2 \over \sigma_f^3} (f(x_0) - f^*)$.
As we saw, the first and the second methods   need  $O(\Delta(x_0))$
and   $\widetilde{O}(\Delta^{1/2}(x_0))$ iterations respectively. Let us show that
for our particular subclass of self-concordant functions
these bounds can be significantly improved.

We consider methods based on the {\em cubic regularization} of
the Newton method (CRNM). Let us define quadratic approximation
of $f$ at point $x \in \E$:
$$
\ba{rcl}
Q(x,y) & = & f(x) + \la \nabla f(x), y - x \ra + \half \la
\nabla^2 f(x)(y-x),y-x\ra.
\ea
$$
Then, from \eqref{eq-Hess}, $| f(y) - Q(x,y) |  \leq  {H_f \over 6} \| y - x \|^3,
\quad y \in \E$,
which justifies the cubic regularized Newton step given a parameter $M>0$:
\beq\label{def-TCN}
\ba{rcl}
x_+=T_M(x) & \Def & \arg\min\limits_{y \in \E} \{ Q(x,y) + {1
\over 6} M \| y - x \|^3 \}.
\ea
\eeq
The cubic regularized Newton method (CRNM) iterating these
steps converges  \cite{nesterov2006cubic} for functions satisfying \eqref{eq-Hess} as $O({1 \over k^2})$, where $k$ is the iteration
counter.

Define the region of quadratic convergence of the CRNM in terms of the function value (see (6.4) in \cite{nesterov2008accelerating}):
$$
\ba{rcl}
\Q_f & = & \left\{ x \in \E: \; f(x) - f^* \leq
{\sigma_f^3 \over 2H_f^2} \refEQ{eq-SS} {1 \over 8M_f^2}
\right\}.
\ea
$$
Let us check how many
iterations we need for entering $\Q_f$ by different
schemes based on the 	step \eqref{def-TCN}. Let our method
have the following convergence rate
\beq\label{eq-SRate}
\ba{rcl}
f(x_k) - f^* & \leq & {c H_f \| x_0 - x^* \|^3 \over k^p}
\; \refLE{eq-strong} \; {c H_f \over k^p} \left({2 \over
\sigma_f} (f(x_0) - f^*) \right)^{3/2}\\
& \refEQ{eq-SS} & {2^{5/2}c M_f \over k^p} (f(x_0) -
f^*)^{3/2},
\ea
\eeq
where $c$ is an absolute constant and $p > 0$. Thus, we
need
$$
\ba{rcl}
O\left(\left[M_f^3(f(x_0)-f^*)^{3/2}\right]^{1/p}\right) &
= & O\left( \Delta^{3 \over 2p}(x_0) \right)
\ea
$$
iterations for entering the region of superlinear
convergence $\Q_f$. For the CRNM we have $p =
2$ (see \cite{nesterov2006cubic}). Thus, it ensures complexity
$O(\Delta^{3/4}(x_0))$. For the accelerated CRNM \cite{nesterov2008accelerating} we have $p=3$. Thus, it needs
$O(\Delta^{1/2}(x_0))$ iterations (which is slightly
better than (\ref{eq-CompPF}) by PFS). However, note that there
exists a powerful tool for accelerating these schemes, the
{\em restarting strategy}.

Let us define $k_p$ as the first integer, for which the
right-hand side of inequality (\ref{eq-SRate}) is smaller
than $\half (f(x_0) - f^*)$:
$$
\ba{rcl}
{2^{5/2}c M_f \over k_p^p} (f(x_0) - f^*)^{3/2} & \leq &
\half (f(x_0) - f^*).
\ea
$$
Clearly, $k_p = O\left( \left[M_f
(f(x_0)-f^*)^{1/2}\right]^{1/p}\right) = O\left( \Delta^{1
\over 2p}(x_0) \right)$. This value can be used in the following multi-stage
scheme.
\beq\label{met-Multi}
\ba{|c|}
\hline \\
\mbox{\bf Multi-stage Acceleration Scheme}\\
\\
\hline \\
\quad \mbox{\BMP At the first stage, we perform $t_1 = \lceil k_p
\rceil$ iterations of our method starting from the point
$y_0 = x_0$ and get the point $y_1$, which is the
starting point for the next stage. 

In general, $k$th stage starts from the point $y_{k-1}$
and its length is  $t_k = \left\lceil {k_p \over
2^{(k-1)/(2p)}} \right\rceil$. Method stops when $y_k \in
\Q_f$.\EMP}
\quad\\
\\
\hline
\ea
\eeq
\BT\label{th-Multi}
The total number of stages $T$ in the optimizations
strategy (\ref{met-Multi}) satisfies inequality $T  \leq  4 + \log_2 \Delta(x_0)$.
The total number $N$ of the lower-level iterations  in this
scheme does not exceed $4 + \log_2 \Delta(x_0)+{2^{1/(2p)} \over 2^{1/(2p)}-1}
k_p$.
\ET
\proof
Let us prove by induction that $f(y_{k}) - f^* \leq
(\half)^k (f(y_0) - f^*)$. For $k=0$ this is true. Assume
that this is also true for some $k \geq 0$. Note that
$t_{k+1}^p \geq (\half)^{k/2} k_p^p$. Therefore,
$$
\ba{rcl}
{f(y_{k+1}) - f^* \over f(y_k) - f^*} & \leq & {2^{5/2}c M_f
\over t_{k+1}^p} (f(y_k) - f^*)^{1/2} \; \leq \; {k_p^p
(f(y_k) - f^*)^{1/2} \over 2  t_{k+1}^p (f(x_0) -
f^*)^{1/2}} \; \\
& \leq & \half \left[ {2^k(f(y_k) - f^*) \over
f(x_0) - f^*} \right]^{1/2} \; \leq \; \half.
\ea
$$
Hence, $T$ satisfies inequality
$\left( \half \right)^{T-1} (f(x_0) - f^*) \geq {1 \over 8
M_f^2}$. Finally,
$$
\ba{rcl}
N & = & \sum\limits_{k=1}^T t_k \; \leq \; T + k_p
\sum\limits_{k=0}^{T-1} \left( \half \right)^{k \over 2p}
\; \leq \; T + k_p \sum\limits_{k=0}^{\infty} \left( \half
\right)^{k \over 2p} \; = \; T + {k_p \over 1 -
\left(\half\right)^{1/(2p)}}.
\ea
$$
\qed

Applying Theorem \ref{th-Multi} to different CRNMs, we get the
following complexity bounds.
\begin{itemize}
\item
{\bf Cubic Regularized Newton Method \cite{nesterov2006cubic}.} For this method $p =
2$. Therefore, the complexity bound of this scheme, as
applied in the framework of multi-stage method
(\ref{met-Multi}) is of the order $O(\Delta^{1/4}(x_0))$.
In fact, this method does not need a restarting strategy.
Thus, Theorem \ref{th-Multi} provides the CRNM
with a better way of estimating its rate of convergence.
\item
{\bf Accelerated Cubic Regularized Newton Method \cite{nesterov2008accelerating}.} For this
method $p=3$. Hence, the complexity bound of the
corresponding multi-stage scheme (\ref{met-Multi}) becomes
$O(\Delta^{1/6}(x_0))$.
\item
{\bf Optimal second-order method \cite{monteiro2013accelerated}.} For this
method $p =3.5$. Therefore, the corresponding complexity
bound is $\tilde{O}(\Delta^{1/7}(x_0))$. However, note
that this method includes an expensive line-search
procedure. Consequently, its practical efficiency should
be worse than the efficiency of the method from the
previous item. Note that the theoretical gap in the
complexity estimates of these methods is negligibly small,
of the order of $O(\Delta^{1/42}(x_0))$.
\end{itemize}
\begin{remark}
Note that the knowledge of $f^*$ is not required for the restarting procedure. Indeed, if we know a lower bound $\tilde{f}$ such that $f^* \geq \tilde{f}$, we obtain that after 
$
\tilde{k}_p = O\left( \left[M_f(f(x_0)-\tilde{f})^{1/2}\right]^{1/p}\right)
$ 
steps of the inner method the  right-hand side of inequality (\ref{eq-SRate}) is smaller
than $\half (f(x_0) - f^*)$ since $\tilde{k}_p \geq k_p$. Then, the overall number of lower-level iterations is of the order of $\tilde{k}_p = O\left( \tilde{\Delta}(x_0)\right)^{\frac{1}{2p}}$, where $\tilde{\Delta}(x_0) = M_f^2(f(x_0)-\tilde{f})$.

\end{remark}

As we can see, the methods considered in this section have
much better complexity bounds for problem
(\ref{prob-SCF}) than the methods based on the
framework of self-concordant functions. A possible
explanation of this phenomena is that these methods use a
more precise model of the objective function, which is
based on two independent inequalities~(\ref{eq-strong})
and (\ref{eq-Hess}) instead of single inequality
(\ref{def-SCF}). Nevertheless, the methods in this section rely on Lipschitz property of the Hessian and, thus, are not applicable for a wide class of self-concordant functions, namely, self-concordant barriers, a typical example being log-barrier $-\ln x$. On the contrary, the DNM and the PFS can be still applied for the minimization of self-concordant barriers.

\section*{Funding and conflicts of interests}
This paper has received funding from the European Research Council (ERC) under the European Union’s Horizon 2020 research and innovation program (grant agreement No 788368). It was also supported by Multidisciplinary Institute in Artificial intelligence MIAI@Grenoble Alpes (ANR-19-P3IA-0003). 
The work by P. Dvurechensky was supported by the Deutsche Forschungsgemeinschaft (DFG, German Research Foundation).
The authors declare that they have no conflict of interest.

\bibliographystyle{plain}
\bibliography{PF_PC_paper_3.bbl}   

\begin{thebibliography}{10}

\bibitem{carderera2020second-order}
Alejandro Carderera and Sebastian Pokutta.
\newblock Second-order conditional gradient sliding.
\newblock {\em arXiv:2002.08907}, 2020.

\bibitem{cartis2011adaptive}
Coralia Cartis, Nicholas~IM Gould, and Philippe~L Toint.
\newblock Adaptive cubic regularisation methods for unconstrained optimization.
  part i: motivation, convergence and numerical results.
\newblock {\em Mathematical Programming}, 127(2):245--295, 2011.

\bibitem{cartis2019universal}
Coralia. Cartis, Nick~I. Gould, and Philippe~L. Toint.
\newblock Universal regularization methods: Varying the power, the smoothness
  and the accuracy.
\newblock {\em SIAM Journal on Optimization}, 29(1):595--615, 2019.

\bibitem{conn2000trust}
Andrew Conn, Nicholas Gould, and Philippe Toint.
\newblock {\em Trust Region Methods}.
\newblock Society for Industrial and Applied Mathematics, 2000.

\bibitem{dinh2013inexact}
Quoc~Tran Dinh, Ion Necoara, Carlo Savorgnan, and Moritz Diehl.
\newblock An inexact perturbed path-following method for lagrangian
  decomposition in large-scale separable convex optimization.
\newblock {\em SIAM Journal on Optimization}, 23(1):95--125, 2013.

\bibitem{doikov2023minimizing}
Nikita Doikov.
\newblock Minimizing quasi-self-concordant functions by gradient regularization
  of newton method.
\newblock {\em arXiv:2308.14742}, 2023.

\bibitem{doikov2021minimizing}
Nikita Doikov and Yurii Nesterov.
\newblock Minimizing uniformly convex functions by cubic regularization of
  newton method.
\newblock {\em Journal of Optimization Theory and Applications},
  189(1):317--339, Apr 2021.

\bibitem{dvurechensky2018global}
Pavel Dvurechensky and Yurii Nesterov.
\newblock Global performance guarantees of second-order methods for
  unconstrained convex minimization.
\newblock CORE Discussion Paper 2018/32, CORE UCL, 2018.

\bibitem{dvurechensky2024hessian}
Pavel Dvurechensky and Mathias Staudigl.
\newblock Hessian barrier algorithms for non-convex conic optimization.
\newblock {\em Mathematical Programming}, 2024.
\newblock (accepted), arXiv:2111.00100.

\bibitem{grapiglia2020tensor}
G.~N. Grapiglia and Yu. Nesterov.
\newblock Tensor methods for minimizing convex functions with h\"older
  continuous higher-order derivatives.
\newblock {\em SIAM Journal on Optimization}, 30(4):2750--2779, 2020.

\bibitem{grapiglia2017regularized}
Geovani~Nunes Grapiglia and Yu~Nesterov.
\newblock Regularized newton methods for minimizing functions with {H}\"older
  continuous hessians.
\newblock {\em SIAM Journal on Optimization}, 27(1):478--506, 2017.

\bibitem{hanzely2022damped}
Slavom\'{\i}r Hanzely, Dmitry Kamzolov, Dmitry Pasechnyuk, Alexander Gasnikov,
  Peter Richtarik, and Martin Takac.
\newblock A damped newton method achieves global $o (1/k^2)$ and local
  quadratic convergence rate.
\newblock In S.~Koyejo, S.~Mohamed, A.~Agarwal, D.~Belgrave, K.~Cho, and A.~Oh,
  editors, {\em Advances in Neural Information Processing Systems}, volume~35,
  pages 25320--25334. Curran Associates, Inc., 2022.

\bibitem{kantorovich1948newton}
L.V. Kantorovich.
\newblock On newton’s method for functional equations.
\newblock {\em Dokl. Akad. Nauk SSSR}, 59(7):1237--1240, 1948.

\bibitem{lee2014proximal}
Jason~D. Lee, Yuekai. Sun, and Michael~A. Saunders.
\newblock Proximal newton-type methods for minimizing composite functions.
\newblock {\em SIAM Journal on Optimization}, 24(3):1420--1443, 2014.

\bibitem{li2017inexact}
Jinchao Li, Martin~S. Andersen, and Lieven Vandenberghe.
\newblock Inexact proximal newton methods for self-concordant functions.
\newblock {\em Mathematical Methods of Operations Research}, 85(1):19--41, Feb
  2017.

\bibitem{liu2020newton}
D.~Liu, V.~Cevher, and Q.~Tran-Dinh.
\newblock A newton frank-wolfe method for constrained self-concordant
  minimization.
\newblock {\em arXiv:2002.07003}, 2020.

\bibitem{liu2020inexact}
Deyi Liu and Quoc Tran-Dinh.
\newblock An inexact interior-point lagrangian decomposition algorithm with
  inexact oracles.
\newblock {\em Journal of Optimization Theory and Applications},
  185(3):903--926, Jun 2020.

\bibitem{monteiro2013accelerated}
R.~Monteiro and B.~Svaiter.
\newblock An accelerated hybrid proximal extragradient method for convex
  optimization and its implications to second-order methods.
\newblock {\em SIAM Journal on Optimization}, 23(2):1092--1125, 2013.

\bibitem{nesterov2008accelerating}
Yu. Nesterov.
\newblock Accelerating the cubic regularization of newton's method on convex
  problems.
\newblock {\em Mathematical Programming}, 112(1):159--181, Mar 2008.

\bibitem{nesterov2016local}
Yu. Nesterov and L.~Tunçel.
\newblock Local superlinear convergence of polynomial-time interior-point
  methods for hyperbolicity cone optimization problems.
\newblock {\em SIAM Journal on Optimization}, 26(1):139--170, 2016.

\bibitem{nesterov2004introduction}
Yurii Nesterov.
\newblock {\em Introductory Lectures on Convex Optimization: a basic course}.
\newblock Kluwer Academic Publishers, Massachusetts, 2004.

\bibitem{nesterov1994interior}
Yurii Nesterov and Arkadii Nemirovskii.
\newblock {\em Interior-point polynomial algorithms in convex programming}.
\newblock SIAM, 1994.

\bibitem{nesterov2006cubic}
Yurii Nesterov and Boris Polyak.
\newblock Cubic regularization of newton method and its global performance.
\newblock {\em Mathematical Programming}, 108(1):177--205, 2006.

\bibitem{ODonoghue2016conic}
Brendan O'Donoghue, Eric Chu, Neal Parikh, and Stephen Boyd.
\newblock Conic optimization via operator splitting and homogeneous self-dual
  embedding.
\newblock {\em Journal of Optimization Theory and Applications},
  169(3):1042--1068, Jun 2016.

\bibitem{polyak2019new}
Boris Polyak and Andrey Tremba.
\newblock New versions of newton method: step-size choice, convergence domain
  and under-determined equations.
\newblock {\em Optimization Methods and Software}, 0(0):1--32, 2019.

\bibitem{rodomanov2016superlinearly}
Anton Rodomanov and Dmitry Kropotov.
\newblock A superlinearly-convergent proximal newton-type method for the
  optimization of finite sums.
\newblock In Maria~Florina Balcan and Kilian~Q. Weinberger, editors, {\em
  Proceedings of The 33rd International Conference on Machine Learning},
  volume~48 of {\em Proceedings of Machine Learning Research}, pages
  2597--2605, New York, New York, USA, 20--22 Jun 2016. PMLR.

\bibitem{song2021unified}
Chaobing Song, Yong Jiang, and Yi~Ma.
\newblock Unified acceleration of high-order algorithms under general hölder
  continuity.
\newblock {\em SIAM Journal on Optimization}, 31(3):1797--1826, 2021.

\bibitem{staudigl2020self-concordant}
Mathias Staudigl, Pavel Dvurechensky, Shimrit Shtern, Kamil Safin, and Petr
  Ostroukhov.
\newblock Self-concordant analysis of {F}rank-{W}olfe algorithms.
\newblock In {\em Proceedings of the 37th International Conference on Machine
  Learning}, volume 119 of {\em Proceedings of Machine Learning Research}.
  PMLR, 2020.
\newblock (accepted), arXiv:2002.04320.

\bibitem{sun2018generalized}
Tianxiao Sun and Quoc Tran-Dinh.
\newblock Generalized self-concordant functions: a recipe for newton-type
  methods.
\newblock {\em Mathematical Programming}, 178(1-2):145--213, 2019.

\bibitem{tran-dinh2018new}
Quoc Tran-Dinh, Liang Ling, and Kim-Chuan Toh.
\newblock A new homotopy proximal variable-metric framework for composite
  convex minimization.
\newblock {\em arXiv:1812.05243}, 2018.

\bibitem{zhang2018communication}
Yuchen Zhang and Lin Xiao.
\newblock {\em Communication-Efficient Distributed Optimization of
  Self-concordant Empirical Loss}, pages 289--341.
\newblock Springer International Publishing, Cham, 2018.

\end{thebibliography}


%
%
%
%
%
%
%
%
%
%

\end{document}